\documentclass[12pt]{amsart}
\usepackage{amssymb}
\usepackage{amscd}		
\usepackage{graphicx}		

\oddsidemargin \evensidemargin
   \title{A Combinatorial Formula for the Character of the
          Diagonal Coinvariants}

   \author{J. Haglund}
   \thanks{Work supported by an NSA grant (J.H.)}

   \author{M. Haiman}
   \thanks{Work supported by NSF Grants DMS-0296203 and DMS-0301072 (M.H.)}

   \author{N. Loehr}
   \thanks{Work supported by an NSF Graduate Fellowship (N.L.)}

   \author{J. B. Remmel}
   \thanks{}

   \author{A. Ulyanov}
   \thanks{}

   \address[J.H., N.L., A.U.]{Dept.\ of Mathematics\\
            University of Pennsylvania\\
            Philadelphia, PA}

   \address[M.H.]{Dept.\ of Mathematics\\
            University of California\\
            Berkeley, CA}

   \address[J.B.R.]{Dept.\ of Mathematics\\
                University of California\\
                San Diego, CA}

   \email[J.H.]{jhaglund@math.upenn.edu}
   \email[M.H.]{mhaiman@math.berkeley.edu}
   \email[N.L.]{nloehr@math.upenn.edu}
   \email[J.B.R.]{remmel@math.ucsd.edu}
   \email[A.U.]{ulyanov@math.upenn.edu}

   \date{October 26, 2003; revised March 2, 2004}

\newtheorem{thm}{Theorem}[subsection]
\newtheorem{lemma}[thm]{Lemma}
\newtheorem{prop}[thm]{Proposition}
\newtheorem{cor}[thm]{Corollary}
\newtheorem{conj}[thm]{Conjecture}

\theoremstyle{definition}
\newtheorem{defn}[thm]{Definition}

\theoremstyle{remark}
\newtheorem{problem}{Problem}
\newtheorem*{example}{Example}
\newtheorem*{remark}{Remark}
\newtheorem*{remarks}{Remarks}

\DeclareMathOperator{\ch}{char}
\DeclareMathOperator{\SYT}{SYT}
\DeclareMathOperator{\SSYT}{SSYT}
\DeclareMathOperator{\SSRT}{SSRT}
\DeclareMathOperator{\dinv}{dinv}
\DeclareMathOperator{\inv}{inv}
\DeclareMathOperator{\core}{core}
\DeclareMathOperator{\quot}{quot}
\DeclareMathOperator{\spin}{sp}
\DeclareMathOperator{\cospin}{csp}
\DeclareMathOperator{\smin}{smin}
\DeclareMathOperator{\smax}{smax}

\DeclareMathOperator{\area}{area}
\DeclareMathOperator{\maj}{maj}
\DeclareMathOperator{\comaj}{comaj}
\DeclareMathOperator{\bounce}{bounce}
\DeclareMathOperator{\mmod}{mod}
\DeclareMathOperator{\ev}{ev}

\newcommand{\Acal}{{\mathcal A}}
\newcommand{\Dcal}{{\mathcal D}}
\newcommand{\Fcal}{{\mathcal F}}
\newcommand{\Hcal}{{\mathcal H}}
\newcommand{\Pcal}{{\mathcal P}}

\newcommand{\Sfrak}{{\mathfrak S}}

\newcommand{\CC}{{\mathbb C}}
\newcommand{\NN}{{\mathbb N}}
\newcommand{\ZZ}{{\mathbb Z}}

\newcommand{\xx}{{\mathbf x}}
\newcommand{\yy}{{\mathbf y}}
\newcommand{\boldmu}{{\boldsymbol \mu}}

\newcommand{\defeq}{\underset{\text{def}}{=}}

\newcommand{\ket}[1]{{| #1 \rangle}}

\newcommand{\qbinom}[2]{\genfrac{[}{]}{0pt}{}{#1}{#2}}

\newlength{\cellsize} 
\newsavebox{\cell}
\sbox{\cell}{\begin{picture}(18,18)
\put(0,0){\line(1,0){18}}
\put(0,0){\line(0,1){18}}
\put(18,0){\line(0,1){18}}
\put(0,18){\line(1,0){18}}
\end{picture}}
\newcommand\cellify[1]{\def\thearg{#1}\def\nothing{}%
\ifx\thearg\nothing
\vrule width0pt height\cellsize depth0pt\else
\hbox to 0pt{\usebox{\cell} \hss}\fi%
\vbox to \cellsize{
\vss
\hbox to \cellsize{\hss$#1$\hss}
\vss}}
\newcommand\tableau[1]{\vcenter{\let\\\cr
\baselineskip -16000pt \lineskiplimit 16000pt \lineskip 0pt
\ialign{&\cellify{##}\cr#1\crcr}}}
\begin{document}

\subjclass[2000]{Primary: 05E10; Secondary: 05A30, 20C30}

\begin{abstract}
Let $R_{n}$ be the ring of coinvariants for the diagonal action of the
symmetric group $S_{n}$.  It is known that the character of $R_{n}$ as
a doubly-graded $S_{n}$ module can be expressed using the Frobenius
characteristic map as $\nabla e_{n}$, where $e_{n}$ is the $n$-th
elementary symmetric function, and $\nabla $ is an operator from the
theory of Macdonald polynomials.

We conjecture a combinatorial formula for $\nabla e_{n}$ and prove
that it has many desirable properties which support our conjecture.  In
particular, we prove that our formula is a symmetric function (which
is not obvious) and that it is Schur positive.  These results make use
of the theory of ribbon tableau generating functions of Lascoux,
Leclerc and Thibon.  We also show that a variety of earlier
conjectures and theorems on $\nabla e_{n}$ are special cases of our
conjecture.

Finally, we extend our conjectures on $\nabla e_{n}$ and several of
the results supporting them to higher powers $\nabla ^{m}e_{n}$.
\end{abstract}

\maketitle

\section{Introduction}
\label{intro}

\subsection{}
\label{s:intro}

Let $R_{n}$ be the ring of coinvariants for the diagonal action of the
symmetric group $S_{n}$ on $\CC ^{n}\oplus \CC ^{n}$.  In other words,
\begin{equation}\label{e:Rn-introduced}
R_{n} = \CC [\xx ,\yy ]/I,
\end{equation}
where $\CC [\xx ,\yy ]=\CC [x_{1},y_{1},\ldots,x_{n},y_{n}]$ is the
ring of polynomial functions on $\CC ^{n}\oplus \CC ^{n}$, the
symmetric group acts ``diagonally'' ({\it i.e.}, permuting the $x$ and
$y$ variables simultaneously), and the ideal $I = ((\xx ,\yy )\cap \CC
[\xx ,\yy ]^{S_{n}})$ is generated by all $S_{n}$-invariant
polynomials without constant term.  The $S_{n}$ action respects the
double grading
\begin{equation}\label{e:Rn-grading}
R_{n} = \bigoplus _{r,s}(R_{n})_{r,s}
\end{equation}
given by the $x$ and $y$ degrees.

A formula for the character of $R_{n}$ as a doubly graded
$S_{n}$ module was conjectured in \cite{GaHa96} and proved in
\cite{Hai02}.  The formula expresses the character in terms of
Macdonald polynomials, as follows.  Let $F$ denote the Frobenius
characteristic: the linear map from $S_{n}$ characters to
symmetric functions that sends the irreducible character $\chi
^{\lambda }$ to the Schur function $s_{\lambda }(z)$.  Encoding the
graded character of $R_{n}$ by means of its {\it Frobenius series}
\begin{equation}\label{e:Rn-Frobenius-series}
\Fcal _{R_{n}}(z;q,t) = \sum _{r,s}q^{r}t^{s}F \ch (R_{n})_{r,s},
\end{equation}
its value is given by the following theorem.
\begin{thm}[\cite{Hai02}]\label{thm:FRn=nabla-en}
Let $\nabla $ be the linear operator defined in terms of the
modified Macdonald symmetric functions $\tilde{H}_{\mu }(z;q,t)$ by
\begin{equation}\label{e:nabla}
\nabla \tilde{H}_{\mu } = t^{n(\mu )}q^{n(\mu ')} \tilde{H}_{\mu },
\end{equation}
where $\mu $ is a partition of $n$, $\mu '$ is its conjugate and
$n(\mu ) = \sum _{i}(i-1)\mu _{i}$.  Then we have
\begin{equation}\label{e:FRn=nabla-en}
\Fcal _{R_{n}}(z;q,t) = \nabla e_{n}(z),
\end{equation}
where $e_{n}$ is the $n$th elementary symmetric function.
\end{thm}

The operator $\nabla $ has been the subject of a series of theorems
and conjectures of a combinatorial nature
\cite{BeGaHaTe99,GaHa96,Hai94,Hai02} (see also \cite{Hai03y} for an
overview).  Specifically, thanks to results of Garsia and Haiman in
\cite{GaHa96}, Theorem~\ref{thm:FRn=nabla-en} implies that the
dimension of $R_{n}$ is given by
\begin{equation}\label{e:(n+1)^(n-1)}
\dim _{\CC }R_{n} = (n+1)^{(n-1)},
\end{equation}
and that of its subspace $R_{n}^{\epsilon }$ of $S_{n}$-antisymmetric
elements by 
\begin{equation}\label{e:Rn-epsilon}
\dim _{\CC }R_{n}^{\epsilon } = C_{n} = \frac{1}{n+1}\binom{2n}{n},
\end{equation}
the $n$-th {\it Catalan number}.  These and other related results
suggest that we should try to understand the rather mysterious
quantity $\nabla e_{n}(z)$ in more combinatorial terms.  Taking a
first step in this direction, Garsia and Haglund
\cite{GaHag01a,GaHag02} gave an explicit combinatorial formula for the
Hall inner product
\begin{equation}\label{e:Cn(q,t)}
C_{n}(q,t) = \langle \nabla e_{n}, e_{n}  \rangle,
\end{equation}
which by Theorem~\ref{thm:FRn=nabla-en} and equation
\eqref{e:Rn-epsilon} is a $q,t$-analog of the Catalan number
$C_{n}(1,1) = C_{n}$.  Building on the Garsia-Haglund formula, Haglund
and Loehr \cite{HagLoe02} conjectured a combinatorial formula for the
Hilbert series of $R_{n}$.  By Theorem~\ref{thm:FRn=nabla-en}, this
Hilbert series is given by
\[
\Hcal _{n}(q,t) = \langle \nabla e_{n}, e_{1}^{n} \rangle = \sum
_{r,s}q^{r}t^{s}\dim (R_{n})_{r,s}.
\]
By \cite{GaHa96}, it was known that $\Hcal _{n}(1,t)$ is a generating
function enumerating {\it parking functions} according to a suitably
defined weight.  The Haglund-Loehr conjecture interprets $\Hcal
_{n}(q,t)$ as a bivariate generating function enumerating parking
functions by the usual weight, together with another statistic
counting certain kinds of inversions (see \S \ref{s:Haglund-Loehr}).

In this paper we conjecture a combinatorial formula for the full
expansion of $\nabla e_{n}(z)$ in terms of monomials, generalizing the
Garsia-Haglund formula for $C_{n}(q,t)$, the Haglund-Loehr conjecture
for $\Hcal _{n}(q,t)$, and a conjecture in \cite{EgHaKiKr03}
expressing $\langle \nabla e_{n}, h_{d}e_{n-d} \rangle$ in terms of
Schr\"oder paths.  We prove that our formula is, as it ought to be, a
symmetric function.  As will be seen, this property of our formula is
not obvious from its definition, but follows from the theory of ribbon
tableau generating functions developed by Lascoux, Leclerc and Thibon
\cite{LaLeTh97,LecThi00}.

By Theorem~\ref{thm:FRn=nabla-en}, $\nabla e_{n}(z)$ is
Schur positive, that is, its coefficients $\langle \nabla
e_{n}(z),s_{\lambda } \rangle$ on the Schur basis belong to $\NN
[q,t]$.  We prove that our conjectured formula is also, as it ought
be, Schur positive.  For this, however, we must rely on an
interpretation of our formula in terms of Kazhdan-Lusztig polynomials,
as in \cite{LecThi00}.  We are unable as yet to provide a
combinatorial interpretation for its Schur function expansion.

Finally, we extend our considerations to higher powers $\nabla
^{m}e_{n}(z)$, giving corresponding conjectured formulas and examining
their properties.

\section{Preliminaries}
\label{preliminaries}

\subsection{$q$-Series notation}
\label{s:q-series}

We use the standard notations:
\begin{align}\label{e:q-things}
(z;q)_{k}&	= (1-z)(1-zq)\cdots (1-zq^{k-1}),\\
[k]_{q} &	= \frac{1-q^{k}}{1-q}, \\
[k]_{q}!&	= (q;q)_{k}/(1-q)^{k} = [k]_{q}\, [k-1]_{q}\, \cdots\,
[1]_{q},\\ 
\qbinom{n}{k}_{q}& = (q^{n-k+1};q)_{k}/(q;q)_{k} =
\frac{[n]_{q}!}{[k]_{q}![n-k]_{q}!},\\
\qbinom{n}{k_{1},\ldots,k_{r}}_{q}& =
\frac{[n]_{q}!}{[k_{1}]_{q}!\cdots  [k_{r}]_{q}!},\quad \text{where
$k_{1}+\cdots +k_{r} = n$}.
\end{align}

\subsection{Partitions and tableaux}
\label{s:partitions-tableaux}

We represent an integer partition as usual by the sequence
\[
\lambda = (\lambda _{1},\ldots,\lambda _{l})
\]
of its parts in decreasing order, and denote its size by
\[
|\lambda | = \sum _{i}\lambda _{i}.
\]
It is understood that $\lambda _{i}=0$ for $i>l$.  We may also write
\[
\lambda = (1^{\alpha _{1}}, 2^{\alpha _{2}}, \ldots )
\]
to indicate the partition with $\alpha _{i}$ parts equal to $i$.
The {\it conjugate partition} $\lambda '$ is defined by 
\[
\lambda '_{i} = \sum _{j\geq i}\alpha _{j}.
\]

The {\it Young diagram} of $\lambda $ is the set $\{(i,j):0\leq
j<\lambda _{i+1} \}\subseteq \NN \times \NN $.  One pictures elements
$(i,j)\in \NN \times \NN $ as boxes or {\it cells}, arranged with the
$i$-axis vertical and the $j$-axis horizontal, so the rows of the
diagram are the parts of $\lambda $.  Abusing notation, we usually
write $\lambda $ both for a partition and its diagram.  A {\it skew}
Young diagram $\lambda /\mu $ is the difference of partition diagrams
$\mu \subseteq \lambda $.  A skew diagram is a {\it horizontal strip}
(resp.\ {\it vertical strip}) if it contains no two cells in the same
column (resp.\ row).

A {\it semistandard Young tableau} of (skew) shape $\lambda $ is a
function $T$ from the diagram of $\lambda $ to the ordered alphabet 
\[
\Acal _{+} = \{1<2<\cdots  \}
\]
which is weakly increasing on each row of $\lambda $ and strictly
increasing on each column.  A semistandard tableau is {\it standard}
if it is a bijection from $\lambda $ to $\{1,2,\ldots,n=|\lambda |
\}$.  More generally, we admit the alphabet
\[
\Acal _{\pm } = \Acal _{+}\cup \Acal _{-} =
\{1<\bar{1}<2<\bar{2}<\cdots \}
\]
of {\it positive} letters $1,2,\ldots$ and {\it negative} letters
$\bar{1},\bar{2},\ldots$.  A {\it super} tableau is a function
$T\colon \lambda \rightarrow \Acal _{\pm }$, weakly increasing on each
row and column, such that the entries equal to $a$ in $T$ occupy a
horizontal strip if $a$ is positive, and a vertical strip if $a$ is
negative.  Thus a semistandard tableau is just a super tableau with
positive entries.  We denote
\begin{align*}
\SSYT (\lambda )& = \{\text{semistandard tableaux $T\colon \lambda
\rightarrow \Acal _{+}$} \}\\
\SSYT _{\pm }(\lambda )& =  \{\text{super tableaux $T\colon \lambda
\rightarrow \Acal _{\pm }$} \}\\
\SSYT (\lambda ,\mu )& = \{\text{semistandard tableaux $T\colon \lambda
\rightarrow \Acal _{+}$ with entries $1^{\mu _{1}}, 2^{\mu _{2}},
\ldots $} \}\\
\SSYT_{\pm } (\lambda ,\mu ,\eta )& = \{\text{super tableaux
$T\colon \lambda \rightarrow \Acal _{\pm }$ with entries $1^{\mu _{1}},
\bar{1}^{\eta _{1}}, 2^{\mu _{2}}, \bar{2}^{\eta _{2}}, \ldots $} \}\\
\SYT (\lambda ) & = \{\text{standard tableaux
$T\colon \lambda \rightarrow \{1,\ldots,n = |\lambda | \}$} \} = \SSYT
(\lambda ,(1^{n})).
\end{align*}

\subsection{Symmetric functions}
\label{s:symmetric}

We follow the notation of \cite{Mac95}, writing $e_{\lambda }$ for the
elementary symmetric functions, $h_{\lambda }$ for the complete
homogeneous symmetric functions, $m_{\lambda }$ for the monomial
symmetric functions, $p_{\lambda }$ for the power-sums and $s_{\lambda
}$ for the Schur functions.  We take these in variables $z =
z_{1},z_{2},\ldots$ so as not to confuse them with the variables $\xx
$, $\yy $ in $R_{n}$.

We write $\langle -,- \rangle$ for the Hall inner product,
defined by either of the identities
\begin{equation}\label{e:hall-inner-product}
\langle h_{\lambda },m_{\mu } \rangle = \delta _{\lambda \mu } =
\langle s_{\lambda },s_{\mu } \rangle.
\end{equation}
We denote by $\omega $ the involution defined by any of the identities
\begin{equation}\label{e:omega}
\omega e_{\lambda } = h_{\lambda };\quad \omega h_{\lambda
}=e_{\lambda };\quad \omega s_{\lambda } = s_{\lambda '}.
\end{equation}

We use square brackets $f[A]$ to denote the plethystic evaluation of a
symmetric function $f$ on a polynomial, rational function or formal
series $A$.  This is defined by writing $f$ in terms of power sums and
then substituting $p_{m}[A]$ for $p_{m}$, where $p_{m}[A]$ is the
result of substituting $a\mapsto a^{m}$ for every indeterminate in
$A$.  The standard $\lambda $-ring identities hold for plethystic
evaluation, {\it e.g.}, $e_{n}[A+B] = \sum _{k}e_{k}[A]e_{n-k}[B]$,
and so forth.  In particular, setting $Z=z_{1}+z_{2}+\cdots $, we have
$f[Z] = f(z)$ for all $f$.  Using this notation, we may write
\begin{equation}\label{e:coproduct-omega}
\omega ^{W}f[Z+W]
\end{equation}
to denote the result of applying $\omega $ to $f[Z+W] =
f(z_{1},z_{2},\ldots, w_{1},w_{2},\ldots)$, considered as a symmetric
function in the $w$ variables with functions of $z$ as coefficients.
Equations \eqref{e:hall-inner-product} and \eqref{e:omega} then imply
that the coefficient of a monomial $z^{\mu }w^{\eta } = z_{1}^{\mu
_{1}}z_{2}^{\mu _{2}}\cdots w_{1}^{\eta _{1}}w_{2}^{\eta _{2}}\cdots $
in $\omega ^{W}f[Z+W]$ is given by
\begin{equation}\label{e:z-mu-w-eta}
\omega ^{W}f[Z+W] \mid _{z^{\mu }w^{\eta }} = \langle f,e_{\eta
}h_{\mu } \rangle.
\end{equation}

If $T$ is a semistandard tableau of (skew) shape $\lambda $, we set
\begin{equation}\label{e:z-power-T}
z^{T} = \prod _{x\in \lambda }z_{T(x)}.
\end{equation}
Then the familiar combinatorial formula for (skew) Schur functions
reads
\begin{equation}\label{e:schur-function}
s_{\lambda }(z) = \sum _{T\in \SSYT(\lambda )}z^{T}.
\end{equation}
Throughout what follows we fix
\begin{equation}\label{e:Z-and-W}
Z = z_{1}+z_{2}+\cdots ,\quad W = w_{1}+w_{2}+\cdots ,
\end{equation}
and make the convention that
\[
\text{\em $z_{\bar{a}}$ stands for $w_{a}$, for every negative letter
$\bar{a}\in \Acal _{\pm }$.}
\]
In particular, this means that if $T\in \SSYT_{\pm }(\lambda ,\mu ,\eta
)$, then $z^{T} = z^{\mu }w^{\eta }$ by definition.  The
``super'' analog of \eqref{e:schur-function} is then
\begin{equation}\label{e:super-schur-function}
\omega ^{W}s_{\lambda }[Z+W] =  \sum _{T\in \SSYT_{\pm }(\lambda )}z^{T},
\end{equation}
which follows immediately from \eqref{e:z-mu-w-eta} and the Pieri rule.

\subsection{Quasisymmetric functions}
\label{s:quasi}

Let $T(x)=a$, $T(y)=a+1$ be consecutive entries in a standard tableau
$T$, with $x=(i,j)$, $y=(i',j')$.  If $j\geq j'$, we say that $a$ is a
{\it descent} of $T$.  The descent set of $T$ is the subset
\[
d(T) = \{a:\text{$a$ is a descent of $T$} \}\subseteq \{1,\ldots,n-1 \}.
\]
Given any subset $D\subseteq \{1,\ldots,n-1 \}$, Gessel's {\it
quasisymmetric function} is defined by the formula
\begin{equation}\label{e:Qn,D}
Q_{n,D}(z) = \sum _{\substack{a_{1}\leq a_{2}\leq \cdots \leq a_{n}\\
a_{i}=a_{i+1}\, \Rightarrow \,  i\not \in D }}
z_{a_{1}}z_{a_{2}}\cdots z_{a_{n}}.
\end{equation}
Here the indices $a_{i}$ belong to the alphabet of positive letters
$\Acal _{+}$.

\begin{prop}[\cite{Ges86}]\label{prop:schur-via-Q}
The (skew) Schur function $s_{\lambda }(z)$ is given in terms of
quasisymmetric functions by
\[
s_{\lambda }(z) = \sum _{T\in \SYT(\lambda )}Q_{|\lambda |,d(T)}(z).
\]
\end{prop}

We will need a ``super'' version of the above proposition.  To this
end, define ``super'' quasisymmetric functions
\begin{equation}\label{e:super-Qn,D}
\tilde{Q}_{n,D}(z,w) = \sum _{\substack{a_{1}\leq a_{2}\leq \cdots \leq a_{n}\\
a_{i}=a_{i+1}\in \Acal _{+}\, \Rightarrow \,  i\not \in D\\
a_{i}=a_{i+1}\in \Acal _{-}\, \Rightarrow \,  i \in D }}
z_{a_{1}}z_{a_{2}}\cdots z_{a_{n}}.
\end{equation}
Here the indices $a_{i}$ range over $\Acal _{\pm }$.  Note that our
convention $z_{\bar{a}} = w_{a}$ remains in force, so the right-hand
side stands for an expression involving both $z$ and $w$ variables.
The next proposition generalizes Proposition~\ref{prop:schur-via-Q}.
We review the well-known proof because later on we shall want to prove
similar results for other kinds of tableaux by appealing to the same
mode of reasoning.

\begin{prop}\label{prop:super-schur-via-Q}
The (skew) super Schur function $\tilde{s}_{\lambda }(z,w) = \omega
^{W}s_{\lambda }[Z+W]$ is given in terms of super quasisymmetric
functions by
\begin{equation}\label{e:super-schur-via-Q}
\tilde{s}_{\lambda }(z,w) = \sum _{T\in \SYT(\lambda
)}\tilde{Q}_{|\lambda |,d(T)}(z,w).
\end{equation}
In particular, Proposition~\ref{prop:schur-via-Q} follows on setting
$w=0$.
\end{prop}

\begin{proof}
If $\nu $ is a horizontal strip, there is a unique labelling of the
cells of $\nu $ to form a standard tableau with no descents (namely,
label the cells in increasing order by columns).  Symmetrically, if
$\nu $ is a vertical strip, there is a unique standard tableau on $\nu
$ with descents at every position.  From these observations it follows
that every super tableau $T$, say of shape $\lambda $, has a unique
{\it standardization} $S$ such that $S$ is standard, $T\circ S^{-1}$
is a weakly increasing function, and if $T\circ S^{-1}(j)=T\circ
S^{-1}(j+1)=\cdots =T\circ S^{-1}(k) = a$, then $\{j,\ldots,k-1 \}\cap
d(S)$ is empty if $a$ is positive, and equal to $\{j,\ldots,k-1 \}$ if
$a$ is negative.

Conversely, the shape of a standard tableau with no descents can only
be a horizontal strip, and symmetrically, the shape of a tableau with
descents at every position can only be a vertical strip.  It follows
that for a given standard tableau $S$ of shape $\lambda $, if
$T'\colon \{1,\ldots,n \} \rightarrow \Acal _{\pm }$ is a weakly
increasing function that satisfies the conditions above, then $T =
T'\circ S$ is a super tableau, and its standardization is equal to
$S$.  Since the sum of $z^{T}$ over all such $T$ is equal to
$\tilde{Q}_{|\lambda |,d(S)}(z,w)$, it follows that
\eqref{e:super-schur-via-Q} is just another way of writing
\eqref{e:super-schur-function}.
\end{proof}

\begin{cor}\label{cor:superization}
Let $f(z)$ be any symmetric function homogeneous of degree $n$,
written in terms of quasisymmetric functions as
\begin{equation}\label{e:f-by-Q}
f(z) = \sum _{D} c_{D}Q_{n,D}(z).
\end{equation}
Then its ``superization'' $\tilde{f}(z,w) = \omega ^{W}f[Z+W]$ is given by
\begin{equation}\label{e:super-f-by-Q}
\tilde{f}(z,w) = \sum _{D} c_{D}\tilde{Q}_{n,D}(z,w).
\end{equation}
\end{cor}

\begin{proof}
The quasisymmetric functions $Q_{n,D}(z)$ are linearly independent;
hence the coefficients $c_{D}$ are uniquely determined by $f$, and the
right-hand side of \eqref{e:super-f-by-Q} depends linearly on $f$.
When $f$ is a Schur function, the result follows from
Propositions~\ref{prop:schur-via-Q} and \ref{prop:super-schur-via-Q}.
This implies the result for all $f$ by linearity.
\end{proof}

\section{The main conjecture}
\label{main}

\subsection{}
\label{s:main}

\begin{figure}
\begin{center}
\includegraphics{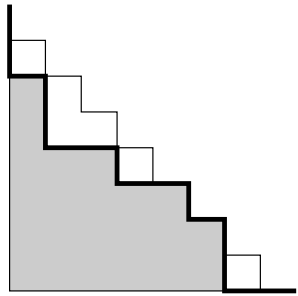}
\caption{\label{fig:dyck}%
A partition $\lambda \subseteq \delta _{n}$ (shaded) and its Dyck path
(heavy line).}
\end{center}
\end{figure}

Fix $n$ and let 
\begin{equation}\label{e:delta-n}
\delta _{n} = (n-1,n-2,\ldots,1,0)
\end{equation}
be the ``staircase'' partition.  Let $\lambda \subseteq \delta _{n}$
be a partition whose diagram is contained in the staircase.  Note that
the outer boundary of $\lambda $, together with segments along the $i$
and $j$-axes, can be identified with a {\it Dyck path}: a lattice path
from $(n,0)$ to $(0,n)$ by steps of the form $(-1,0)$ (south) and
$(0,1)$ (east) that never goes above the diagonal line $i+j=n$.
Figure~\ref{fig:dyck} illustrates this.  The number of Dyck paths, or
of partitions $\lambda \subseteq \delta _{n}$, is the Catalan number
$C_{n}$.

\begin{figure}
\begin{center}
\includegraphics{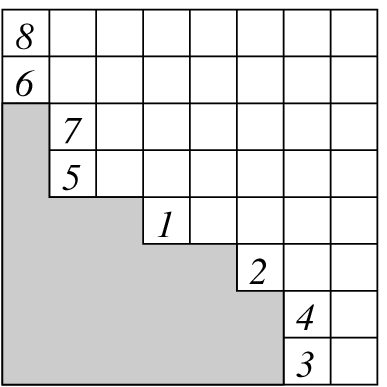}
\caption{\label{fig:park}%
A standard tableau of shape $(\lambda +(1_{n}))/\lambda $.}
\end{center}
\end{figure}

Let $T$ be a semistandard tableau of skew shape $(\lambda
+(1^{n}))/\lambda $, that is, the vertical strip formed by the cells
$(i,\lambda _{i+1})$ for $i=0,1,\ldots,n-1$.  For every cell $x=(i,j)\in
\NN \times \NN $ let $d(x) = i+j$, so $d(x)=k$ means that $x$ is on
the $k$-th diagonal.  Given two entries $T(x)=a$ and $T(y)=b$ of $T$
with $a<b$, put $x=(i,j)$, $y=(i',j')$.  We say that these two entries
form a {\it d-inversion} if either
\begin{itemize}
\item [(i)] $d(y)=d(x)$ and $j>j'$, or
\item [(ii)] $d(y)=d(x)+1$ and $j<j'$.
\end{itemize}
Set
\begin{equation}\label{e:dinv}
\dinv (T) = \text{number of d-inversions in $T$}.
\end{equation}
For example, the tableau $T$ in Figure~\ref{fig:park}
 has $\dinv (T) = 8$, with d-inversions
formed by the pairs of entries
$(8,2)$, $(8,4)$, $(6,1)$, $(6,3)$, $(7,1)$, $(7,3)$, $(4,1)$, and $(2,1)$.

\begin{defn}\label{defn:Dn}
\begin{equation}\label{e:Dn}
D_{n}(z;q,t) = \sum _{\lambda \subseteq \delta _{n}}\; \sum _{T\in \SSYT
(\lambda +(1^{n})/\lambda )} t^{|\delta_{n} / \lambda|} q^{\dinv
(T)}z^{T}.
\end{equation}
\end{defn}

\begin{conj}\label{conj:main}
We have the identity
\begin{equation}\label{e:main}
\nabla e_{n}(z) = D_{n}(z;q,t).
\end{equation}
Equivalently, for all $\mu $ we have 
\begin{equation}\label{e:main-equiv}
\langle \nabla e_{n}, h_{\mu } \rangle = \sum _{\lambda \subseteq
\delta _{n}}\; \sum _{T\in \SSYT (\lambda +(1^{n})/\lambda ,\, \mu )}
t^{|\delta_{n} / \lambda|} q^{\dinv (T)}.
\end{equation}
\end{conj}

The bulk of our work in this paper will serve to show that
Conjecture~\ref{conj:main} is consistent with previously known or
conjectured properties of $\nabla e_{n}$.  The most basic such
property is that $\nabla e_{n}$ is a symmetric function.

\begin{thm}\label{thm:symmetry-positivity}
The quantity $D_{n}(z;q,t)$ is a symmetric function in $z$, and it is
Schur positive, {\em i.e.}, $\langle D_{n}(z;q,t),s_{\mu }(z)
\rangle\in \NN [q,t]$ for all $\mu $.  In fact, each term
\begin{equation}\label{e:lambda-term}
D_{n}^{\lambda }(z;q) = \sum _{T\in \SSYT (\lambda +(1^{n})/\lambda
)} q^{\dinv (T)}z^{T}
\end{equation}
is individually symmetric and Schur positive.
\end{thm}

Note that it is not at all obvious from the definition that
$D_{n}(z;q,t)$ is symmetric.  Its symmetry is equivalent to the
assertion that the right-hand side of \eqref{e:main-equiv} does not
depend on the order of the parts of $\mu $.  Our proof uses Lascoux,
Leclerc and Thibon's theory of spin generating functions for ribbon
tableaux.  We will give a synopsis of their theory and prove
Theorem~\ref{thm:symmetry-positivity} in \S \ref{llt}.  For now, we
take the theorem for granted, and explore what can be deduced by more
elementary means.

\subsection{Superization}
\label{s:superization}

Our first goal is to show that Conjecture~\ref{conj:main} implies a
seemingly stronger formula, giving the superization $\omega ^{W}(\nabla
e_{n})[Z+W]$.  We  extend the definition of $\dinv (T)$ to super
tableaux $T$ as follows.  Let $x,y\in (\lambda +(1^{n}))/\lambda $
satisfy condition (i) or (ii) in the definition of d-inversion,
above.  Then entries $T(x)=a$, $T(y)=b$ with $a,b\in \Acal _{\pm }$
form a d-inversion in $T$ if $a<b$ or if $a=b\in \Acal _{-}$ is a
negative letter.

\begin{thm}\label{thm:super-Dn}
The superization $\tilde{D}_{n}(z,w;q,t) = \omega ^{W}D_{n}[Z+W;q,t]$
is given by
\begin{equation}\label{e:super-Dn}
\tilde{D}_{n}(z,w;q,t) = \sum _{\lambda \subseteq \delta _{n}}\; \sum
_{T\in \SSYT_{\pm } (\lambda +(1^{n})/\lambda )} t^{|\delta_{n} /
\lambda|} q^{\dinv (T)}z^{T}.
\end{equation}
Equivalently, Conjecture~\ref{conj:main} implies
\begin{equation}\label{e:super-main-equiv}
\langle \nabla e_{n}, e_{\eta }h_{\mu } \rangle = \sum _{\lambda
\subseteq \delta _{n}}\; \sum _{T\in \SSYT_{\pm } (\lambda
+(1^{n})/\lambda ,\, \mu, \,\eta )} t^{|\delta_{n} / \lambda|}
q^{\dinv (T)}.
\end{equation}
\end{thm}

\begin{proof}
Consider the total ordering of $\NN \times \NN $ ({\it reverse
diagonal lexicographic order}) defined by
\[
x <_{d} y \quad \text{if $d(x)>d(y)$, or $d(x)=d(y)$ and $j<j'$,
where $x=(i,j)$, $y=(i',j')$}.
\]
On every (skew) shape $\nu $ there is a unique standard tableau whose
labels are decreasing with respect to $<_{d}$; and there is a tableau
with increasing labels, also unique, if and only if all cells of $\nu
$ are in distinct rows and columns.  Suppose now that $\nu $ is
contained in $(\lambda +(1^{n}))/\lambda $.  Then $\nu $ is a vertical
strip {\it a fortiori}, and it is a horizontal strip if and only if
its cells are in distinct rows and columns.  Hence on $\nu $ there
exists a $<_{d}$-decreasing (resp.\ $<_{d}$-increasing) standard
tableau if and only if $\nu $ is a vertical (resp.\ horizontal) strip.
In either case the tableau in question is unique.

Define $a$ to be a d-descent of a standard tableau $S \in \SYT
(\lambda +(1^{n})/\lambda )$ if $S(x)=a$, $S(y)=a+1$ with $x >_{d} y$,
and denote by $dd(S)$ the set of d-descents of $S$.  Define the
standardization of a super tableau $T \in \SSYT (\lambda
+(1^{n})/\lambda )$ to be the unique standard tableau $S$ such that
$T\circ S^{-1}$ is weakly increasing, and if $T\circ S^{-1}(j)=T\circ
S^{-1}(j+1)=\cdots =T\circ S^{-1}(k) = a$, then $\{j,\ldots,k-1 \}\cap
dd(S)$ is empty if $a$ is positive, and equal to $\{j,\ldots,k-1 \}$
if $a$ is negative.  The proof of
Proposition~\ref{prop:super-schur-via-Q} again goes through to show
that the standardization in this new sense exists, and that the sum
$\sum _{T}z^{T}$ over all super tableaux $T$ with standardization $S$
is equal to the super quasisymmetric function
$\tilde{Q}_{n,dd(S)}(z,w)$.

Note that if cells $x,y$ satisfy condition (i) or (ii) in the
definition of d-inversion, then $x >_{d} y$.  In particular, a
standard tableau labelled in $<_{d}$-increasing order has no
d-inversions, while one labelled in $<_{d}$-decreasing order has a
d-inversion in every such pair of cells $x,y$.  With this in mind, we
see that if $T$ is a super tableau and $S$ its standardization, then
$\dinv (S) = \dinv (T)$.  This yields the formula
\begin{equation}\label{e:super-Dn-by-Q}
\sum _{T\in \SSYT_{\pm } (\lambda +(1^{n})/\lambda )} q^{\dinv
(T)}z^{T} = \sum _{S\in \SYT (\lambda +(1^{n})/\lambda )}
q^{\dinv (S)} \tilde{Q}_{n,dd(S)}(z,w).
\end{equation}
Setting $w=0$, we obtain the quasisymmetric function expansion
\begin{equation}\label{e:Dn-by-Q}
D_{n}(z;q,t) = \sum _{\lambda \subseteq \delta _{n}}\; \sum _{S\in
\SYT (\lambda +(1^{n})/\lambda )} t^{|\delta _{n}/\lambda |} q^{\dinv
(S)} Q_{n,dd(S)}(z).
\end{equation}
By \eqref{e:super-Dn-by-Q}, the right-hand side of
\eqref{e:super-Dn} is the superization of this, and the
theorem now follows from Theorem~\ref{thm:symmetry-positivity} and
Corollary~\ref{cor:superization}.
\end{proof}

\subsection{Shuffle formulation}
\label{s:shuffle}
Recall that a {\it parking function} on $n$ cars is a function
$f\colon \{1,\ldots,n \}\rightarrow \{1,\ldots,n \}$ satisfying
\begin{equation}\label{e:parking-function}
|f^{-1}(\{1,\ldots,k \})|\geq k\quad \text{for all $k=1,\ldots,n$}.
\end{equation}
For every function $f\colon \{1,\ldots,n \}\rightarrow \NN _{>0}$
there is a unique partition $\lambda $ with at most $n$ parts and a
standard tableau $T$ of shape $(\lambda +(1^{n}))/\lambda $ such that
each entry $a$ of $T$ lies in column $f(a)-1$.  Namely, the parts of
$\lambda $ are the values $f(i)-1$, and the entries of $T$ in the
cells $(i-1,\lambda _{i})$ for which $\lambda _{i} = j-1$ are the
elements of $f^{-1}(\{j \})$.  It is easy to see that $f$ is a parking
function if and only $\lambda \subseteq \delta _{n}$.

Let $f$ be a parking function encoded by $\lambda \subseteq \delta
_{n}$ and $T\in \SYT (\lambda +(1^{n})/\lambda )$.  Reading off the
entries of $T$ in $<_{d}$-increasing order yields a permutation
$w(f)$, with descent set
\begin{equation}\label{e:D(w)=dd(T)}
D(w(f)^{-1}) = dd(T).
\end{equation}
For example, for the parking function encoded by the tableau in
Figure~\ref{fig:park}, we have $w(f)=82467135$ and $D(w(f)^{-1}) =
dd(T) = \{1,3,5,7 \}$.

Say that a permutation $w$ is a $\mu, \eta $-{\it shuffle} if its
inverse is the concatenation of alternately increasing and decreasing
sequences of lengths $\mu _{1}, \eta _{1}, \mu _{2},\eta _{2},
\ldots$.  Define the {\it area} $a(f)$ to be $|\delta _{n}/\lambda | =
\binom{n+1}{2}-\sum _{i}f(i)$; this is equal to the traditional weight
of the parking function, as in \cite{GaHa96,HagLoe02,Hai94}.  Define
$\dinv (f) = \dinv (T)$.  Then we have the following corollary to
Theorem~\ref{thm:symmetry-positivity} and the proof of
Theorem~\ref{thm:super-Dn}.

\begin{cor}\label{cor:shuffle}
Conjecture~\ref{conj:main} implies that  $\langle \nabla
e_{n},e_{\eta }h_{\mu } \rangle$ is the generating function
\begin{equation}\label{e:shuffle}
\sum _{f}t^{a(f)}q^{\dinv (f)},
\end{equation}
summed over parking functions $f$ such that $w(f)$ is a $\mu ,\eta
$-shuffle.  Even without assuming Conjecture~\ref{conj:main}, the
above sum is independent of the order of the parts of $\eta $ and $\mu
$.
\end{cor}

\begin{remark}
The sum in \eqref{e:shuffle} is also independent of the way in which
$\mu $ and $\eta $ are interleaved, as can be seen by setting some of
the parts in the standard interleaving $\mu _{1},\eta _{1},\mu
_{2},\eta _{2},\ldots$ to zero.
\end{remark}

\section{Specializations} \label{specializations}

\subsection{Value at $q=1$}
\label{s:q=1}

By \cite[Thms.~2.1, 2.2 and 3.6]{GaHa96}, we have the formula
\begin{equation}\label{e:nabla-en-q=1}
\nabla e_{n}(z)|_{q=1} = \sum _{\lambda \subseteq \delta _{n}}
t^{|\delta _{n}/\lambda |} e_{\alpha }(z),
\end{equation}
where $\lambda =(0^{\alpha _{0}},1^{\alpha _{1}},2^{\alpha _{2}},\ldots
)$ with $\alpha _{0}$ defined to make $\sum _{i}\alpha _{i} = n$.
We verify that Conjecture~\ref{conj:main} is consistent with this formula.
 
\begin{prop}\label{prop:Dn-q=1}
We have $D_{n}(z;1,t) = \nabla e_{n}(z)|_{q=1}$.
\end{prop}

\begin{proof}
Since the skew Schur function $s_{(\lambda +(1^{n}))/\lambda }$ is
equal to $e_{\alpha }$, this is obvious from \eqref{e:nabla-en-q=1}
and Definition~\ref{defn:Dn}.
\end{proof}

\begin{remark}
Although $D_{n}(z;1,t)$ is trivial to evaluate, $D_{n}(z;q,1)$ is
not---see \S \ref{problems}, Problem~\ref{prob:t=1}.
\end{remark}

\subsection{Value at $t=0$}
\label{s:t=0}

\begin{lemma}\label{lemma:t=0}
We have
\begin{equation}\label{e:t=0}
\nabla e_{n}(z)|_{t=0} = (q;q)_{n}h_{n}[Z/(1-q)].
\end{equation}
\end{lemma}

\begin{proof}
Equation \eqref{e:FRn=nabla-en} implies that $\nabla e_{n}(z)|_{t=0}$
is the Frobenius series in one parameter $q$ of the classical
coinvariant ring $R_{n}\cap \CC [\xx ]$.  By a result of Stanley
\cite{Stanley79}, this is given by the right-hand side of \eqref{e:t=0}.
\end{proof}

\begin{prop}\label{prop:t=0}
We have
\begin{equation}\label{e:Dn-t=0}
D_{n}(z;q,0) = \nabla e_{n}(z)|_{t=0}.
\end{equation}
\end{prop}

\begin{proof}
Only the term $\lambda =\delta _{n}$ contributes when $t=0$.  Then
every cell $x\in (\lambda +(1^{n}))/\lambda $ is on the diagonal
$d(x)=n$, a tableau $T\in \SSYT (\lambda +(1^{n})/\lambda )$ is just a
word in the alphabet $\Acal _{+}$, and $\dinv (T)$ is its number of
inversions in the ordinary sense.  Hence
\begin{equation}\label{e:Dn-t=0-is-what}
D_{n}(z;q,0) = \sum _{\mu } \qbinom{n}{\mu _{1},\ldots,\mu _{l}}_{q}
m_{\mu }(z).
\end{equation}
By the Cauchy formula,
\begin{equation}\label{e:Cauchy-1}
(q;q)_{n}h_{n}[Z/(1-q)] = \sum _{\mu }(q;q)_{n} h_{\mu }[1/(1-q)]
m_{\mu }(z),
\end{equation}
which is equal to the right-hand side of \eqref{e:Dn-t=0-is-what}.
\end{proof}

\subsection{Value at $q=0$}
\label{s:q=0}

We begin by observing that if $\lambda \subseteq \delta _{n}$ and
$\lambda '$ has distinct parts, then for a parking function $f$
encoded by a tableau $T\in \SYT (\lambda +(1^{n})/\lambda )$, the
permutation $w(f)$ is obtained simply by reading $T$ from the top row
to the bottom.  Fix $w_{0}\in S_{n}$ to be the permutation $w_{0}(i) =
n+1-i$, so $w(f)w_{0}$ is $w(f)$ read backwards.

\begin{lemma}\label{lemma:dinv=0}
A standard tableau $T\in \SYT (\lambda +(1^{n})/\lambda )$ has $\dinv
(T) = 0$ if and only if $\lambda '$ has distinct parts and the descent
set of $w(f)w_{0}$ is the set of parts of $\lambda '$, where $f$ is
the parking function encoded by $T$.
\end{lemma}

\begin{proof}
Suppose $\lambda '$ does not have distinct parts.  Then the associated
Dyck path contains two or more consecutive horizontal steps, not on
the $j$-axis, and there is a cell $x\in (\lambda +(1^{n}))/\lambda $
bordering the vertical step which follows them.  There is at least one
other cell $y\in (\lambda +(1^{n}))/\lambda $ with $d(y) = d(x)$ and
$y$ to the left of $x$; fix $y$ to be the rightmost of these.  By
assumption, $y$ and $x$ are not consecutive on their common diagonal.
Hence the cell $z$ directly below $y$ in $\NN \times \NN $ is also in
$(\lambda +(1^{n}))/\lambda $; otherwise $y$ would not have been the
rightmost cell.  No matter what the entries of $T$ in the three cells
$x$, $y$, $z$ are, they form at least one d-inversion.

We have shown that $\dinv (T) = 0$ implies that $\lambda '$ has
distinct parts.  Once this is given, it is easy to see that $\dinv (T)
= 0$ if and only if, in addition, the descents of $w(f)w_{0}$ are the
parts of $\lambda '$.
\end{proof}

\begin{prop}\label{prop:q=0}
We have
\begin{equation}\label{e:Dn-q=0}
D_{n}(z;0,t) = \nabla e_{n}(z)|_{q=0}.
\end{equation}
\end{prop}

\begin{proof}
Recall that the major index $\maj (w)$ of a permutation or a tableau
is defined as the sum of its descents.  Every permutation $w$ occurs
uniquely as $w(f)$ for a shape $\lambda $, tableau $T$ and parking
function $f$ satisfying the conditions in Lemma~\ref{lemma:dinv=0}.
Moreover, we have $|\delta _{n}/\lambda | = \maj (w_{0}w(f)w_{0})$.
Lemma~\ref{lemma:dinv=0} and equations \eqref{e:Dn-by-Q} and
\eqref{e:D(w)=dd(T)} therefore imply that
\begin{equation}\label{e:q=0}
D_{n}(z;0,t) = \sum _{w\in S_{n}} t^{\maj (w_{0}ww_{0})}
Q_{n,D(w^{-1})}(z). 
\end{equation}
Recall that if $(P_{w},Q_{w})$ is the pair of standard tableaux
associated to $w$ by the Schensted correspondence, then $d(Q_{w}) =
D(w)$ and $d(P_{w}) = D(w^{-1})$.  Recall further that
$Q_{w_{0}ww_{0}} = \ev (Q_{w})$, where $\ev $ is the {\it evacuation}
operator of Sch\"utzenberger \cite{Sch77}.  Hence we can rewrite
\eqref{e:q=0} as
\begin{equation}\label{e:schensted}
D_{n}(z;0,t) = \sum _{(P_{w},Q_{w})} t^{\maj (\ev Q_{w})} Q_{n,d(P_{w})}(z).
\end{equation}
Recall \cite[Prop.~7.19.11]{Sta99} that
\begin{equation}\label{e:schur-special}
(t;t)_{n}s_{\lambda }(1/(1-t)) = \sum _{T\in \SYT (\lambda )} t^{\maj (T)}. 
\end{equation}
Using this and Proposition~\ref{prop:schur-via-Q}, we see that the
right-hand side of \eqref{e:schensted} decomposes as
\begin{equation}\label{e:Dn-q=0-final}
\sum _{\lambda } (t;t)_{n}s_{\lambda }(1/(1-t)) s_{\lambda }(z),
\end{equation}
which is equal to $\nabla e_{n}(z)|_{q=0}$ by Lemma~\ref{lemma:t=0}
and the Cauchy formula.
\end{proof}

\subsection{The $q,t$-Catalan formula}
\label{s:qt-Catalan}

Haglund conjectured \cite{Haglund03a} and with Garsia proved
\cite{GaHag01a,GaHag02} a formula for the $q,t$-Catalan polynomial
\begin{equation}\label{e:qt-Catalan}
C_{n}(q,t) = \langle \nabla e_{n}, e_{n} \rangle .
\end{equation}
Haglund and Loehr \cite{HagLoe02} later showed that the formula of
\cite{GaHag01a,GaHag02} can also be written in the form
\begin{equation}\label{e:qt-Catalan-theorem}
C_{n}(q,t) = \sum _{\lambda \subseteq \delta _{n}} t^{|\delta
_{n}/\lambda |}q^{\dinv (\lambda )},
\end{equation}
where $\dinv (\lambda )$ is defined to be $\dinv (T)$ for the super
tableau of shape $(\lambda +(1^{n}))/\lambda$ whose every entry is
$\bar{1}$.  It is immediate from \eqref{e:super-main-equiv} that this
coincides with $\langle D_{n}(z;q,t),e_{n} \rangle$.

\begin{figure}
\begin{center}
\includegraphics{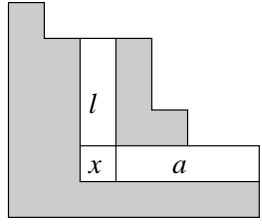}
\caption{\label{fig:arm-leg}%
The arm and leg of a cell $x$ in a Young diagram.}
\end{center}
\end{figure}

The formulation of \eqref{e:qt-Catalan-theorem} in terms of $\dinv
(\lambda )$ was motivated by a conjecture of Haiman, whose original
form was slightly different.  Namely,
\begin{equation}\label{e:qt-Catalan-hook}
C_{n}(q,t) = \sum _{\lambda \subseteq \delta _{n}} t^{|\delta
_{n}/\lambda |}q^{b(\lambda )},
\end{equation}
where $b(\lambda )$ is the number of cells $x\in \lambda $ such that
\begin{equation}\label{e:arm-leg}
l(x)\leq a(x)\leq  l(x)+1.
\end{equation}
Here the {\it arm} $a(x)$ (resp.\ {\it leg} $l(x)$) is the number of
cells in the hook of $x$ that are in the same row (resp.\ column) as
$x$, excluding $x$ itself---see Figure~\ref{fig:arm-leg}.  To tie
\eqref{e:qt-Catalan-theorem} and \eqref{e:qt-Catalan-hook} together,
let us show that in fact $b(\lambda ) = \dinv (\lambda )$.

\begin{lemma}\label{lemma:b-lambda}
The number $b(\lambda )$ defined above is equal to the number of pairs
of cells $x,y\in (\lambda +(1^{n}))/\lambda$ satisfying condition (i)
or (ii) in the definition of d-inversion in \S \ref{s:main}.
\end{lemma}

\begin{proof}
Let $x$ be a cell of $\lambda $, let $u$ be the cell of $(\lambda
+(1^{n}))/\lambda $ just outside the end of the arm of $x$, and let
$t$ be the cell just outside the end of the leg of $x$.  Travelling
along the diagonal $i+j = d(t)$, starting at $t$ and moving in the
increasing $i$ direction, let $v$ be the first cell of $(\lambda
+(1^{n}))/\lambda$ encountered (it always exists).  Then $a(x) = l(x)$
if and only if $u$, $v$ satisfy (i) in the definition of d-inversion,
while $a(x)=l(x)+1$ if and only if $v$, $u$ satisfy (ii).  Moreover,
every pair of cells satisfying (i) or (ii) arises uniquely in this
way.  A fully detailed argument in a more general setting will be
given in the proof of Lemma~\ref{lemma:bm-lambda}, below.
\end{proof}

\subsection{The Haglund-Loehr conjecture}
\label{s:Haglund-Loehr}

The Hilbert series of $R_{n}$ is given by $\Hcal _{n}(q,t) = \langle
\nabla e_{n}, e_{1}^{n} \rangle$.  Conjecture~\ref{conj:main} implies
that this is equal to
\begin{equation}\label{e:Dn-hilbert}
\langle D_{n}(z;q,t),e_{1}^{n} \rangle = \sum _{\lambda \subseteq
\delta _{n}} \;\sum _{T\in \SYT (\lambda +(1^{n})/\lambda )} t^{|\delta
_{n}/\lambda |} q^{\dinv (T)}.
\end{equation}
If desired, one may express the same thing as a sum $\sum
_{f}t^{a(f)}q^{\dinv (f)}$ over all parking functions $f$ on $n$ cars.
It is none other than the value given for $\Hcal _{n}(q,t)$ by a
conjecture of Haglund and Loehr in \cite{HagLoe02}.  Thus the
Haglund-Loehr conjecture is an immediate consequence of
Conjecture~\ref{conj:main}.

\subsection{Fermionic formula}
\label{s:fermionic}

The original ``fermionic formula'' is that of Kerov, Kirillov and
Reshetikhin \cite{KerKirResh86,KirResh86} giving the $q$-Kostka
coefficient $K_{\lambda \mu }(q)$ as a sum of products of $q$-binomial
coefficients.  By analogy, we use the same terminology for an
expansion of a $q,t$-quantity as a sum of powers of $t$ times products
of $q$-binomial coefficients.  Haglund \cite{Haglund03a} gave a
fermionic formula in this sense for $C_{n}(q,t)$.  Here we give a
fermionic formula for $\langle D_{n}(z;q,t), e_{\eta }h_{\mu }
\rangle$.

Let $\sigma = \sigma _1\ldots \sigma _n \in S_n$ be a permutation with
say $k-1$ descents, at positions $r_1<\cdots < r_{k-1}$, and set
$r_{0}=0$, $r_{k}= n$.  Let $A_j=A_j(\sigma) = \{ \sigma _{l}:
r_{j-1}+1\le l \le r_j\}$ be the $j$th ``run" of $\sigma$.  Denote by
$F(\sigma)$ the set of all parking functions $f$ whose encoding
tableau $T$ has the property that $A_{j}$ is the set of entries of $T$
in cells $x$ on the diagonal $d(x) = n+1-j$.
Define
\begin{equation}\label{62}
H(\sigma ;q,t) = \sum _{f \in F(\sigma)}
q^{\dinv (f)}t^{a(f)}.
\end{equation}
With this notation,  \cite[Thm.~1]{HagLoe02} can
be formulated as follows:
\begin{equation}\label{63}
H(\sigma;q,t) = t^{\comaj (\sigma)}\prod _{i=2}^{n}
[v(\sigma ,i)+\chi (i\le r_1)]_{q}, 
\end{equation}
where $\chi(a\leq b)$ equals $1$ if $a\leq b$ and $0$ otherwise,
$\comaj (\sigma) =\sum _{i=1}^{k-1}(n-r_{i}) = \maj (\sigma w_{0})$,
and $v(\sigma ,i)$ is the largest value of $p$, $1\le p \le i-1$, for
which the sequence
\begin{equation}
(\sigma _{i-p} - \sigma _i \mmod   n),\; 
(\sigma _{i-p+1} - \sigma _i \mmod n), \;\ldots, \;
(\sigma _{i-1} -\sigma _i \mmod n)
\end{equation}
is increasing (in other words, $\sigma _{i-p},\ldots,\sigma _{i}$ is a
rotation of an increasing sequence).
Hence $\langle D_{n}(z;q,t), e_{1}^{n} \rangle$, {\it i.e.}, our
conjectured value for the Hilbert series $\Hcal _{n}(q,t) = \langle
\nabla e_{n},e_{1}^{n} \rangle$ of $R_{n}$, is the sum of the
right-hand side of \eqref{63} over all $\sigma \in S_n$.

Let
\begin{equation}\label{64}
H^{\mu}(\sigma; q,t) = \sum _{f \in F^{\mu}(\sigma)}
q^{\dinv (f)}t^{a(f)},
\end{equation}
where $F^{\mu}(\sigma)$ is the set of all parking functions $f\in
F(\sigma)$ such that $w(f)$ is a $\mu, \emptyset$-shuffle.  Note that
if $f\in F(\sigma )$, then $\sigma $ is the permutation obtained by
sorting each block of $w(f)$ contributed by the cells on one diagonal
$d(x)=n+1-j$ in the tableau $T$ that encodes $f$ (these blocks are the
sets $A_{j}$).  Clearly $H^{\mu}(\sigma ;q,t)=0$ if $\sigma$ is not a
$\mu, \emptyset$-shuffle.  Otherwise, if $\sigma $ is a $\mu,
\emptyset$-shuffle, define $B_{j} = \{M_{j-1}+1,\ldots,M_{j} \}$,
where $M_{j} = \mu _{1}+\cdots +\mu _{j}$.  In other words, the
collection $\{B_{j}:1\leq j\leq l(\mu ) \}$ is the partition of
$\{1,\ldots,n \}$ into blocks of consecutive integers of sizes $\mu
_{j}$.

Note that if $f \in F^{\mu}(\sigma)$, then the elements of $A_{i}\cap
B_{j}$ occur in $w(f)$ in increasing order for each $i$ and $j$.  Let
$\Sfrak \subseteq S_{n}$ be the subgroup consisting of permutations
$\tau $ that map each set $A_{i}\cap B_{j}$ into itself.  Then, given
any $g \in F(\sigma)$, there is a unique $f\in F^{\mu }(\sigma )$ and
$\tau \in \Sfrak $ such that $w(g) = \tau \circ w(f)$.
Hence,
\begin{equation}\label{e:obvious-q=1}
(\prod _{i,j} b_{i,j}!) H^{\mu}(\sigma; 1,t)= H(\sigma ;1,t),
\end{equation}
where $b_{i,j} = |A_{i}\cap B_{j}|$.  Moreover, taking into account
the definition of $\dinv $, it is clear that
\begin{equation}\label{e:obvious-q=q}
 H^{\mu}(\sigma; q,t)= \frac{H(\sigma
;q,t)}{\prod _{i,j} [b_{i,j}]_{q}!}.
\end{equation}

Now set $V_{i,j} = v(\sigma ,k)+\chi (k\leq r_{1})$, where $\sigma
_{k}$ is the largest element of $A_{i}\cap B_{j}$.  Using the fact
that the elements of $A_{i}\cap B_{j}$ form an increasing sequence of
adjacent, consecutive integers in $\sigma $, the definition of
$v(\sigma ,i)$ implies that 
\begin{equation}\label{e:H-reformulated}
 \frac{H(\sigma ;q,t)}{\prod _{i,j} [b_{i,j}]_{q}!} = t^{\comaj
(\sigma)}\prod _{i,j} \qbinom{V_{i,j}}{b_{i,j}}_{q}.
\end{equation}
Combining this with \eqref{e:obvious-q=q} yields the fermionic formula
\begin{equation}\label{e:fermionic}
\langle D_{n}(z;q,t),h_{\mu }  \rangle = \sum _{\substack{\sigma \in S_{n}\\
\text{$\sigma $ is a $\mu, \emptyset $-shuffle}}} t^{\comaj
(\sigma)}\prod _{i,j} \qbinom{V_{i,j}}{b_{i,j}}_{q}.
\end{equation}

More generally, there is a similar formula for $\langle
D_{n}(z;q,t),e_{\eta }h_{\mu } \rangle$.  Set $M_{j} = \mu _{1}+\cdots
+ \mu _{j}$ as before, and $E_{j}= \eta _{1}+\cdots +\eta _{j}$.  In
this setting we redefine $B_{j} =
\{M_{j-1}+E_{j-1}+1,\ldots,M_{j}+E_{j-1} \}$ and set $C_{j} =
\{M_{j}+E_{j-1}+1,\ldots,M_{j}+E_{j} \}$.  We define $A_{j}$ to be
the $j$-th run of $\sigma $, just as before.
Let $\tilde{\sigma }$ be
the permutation obtained by reversing each block $A_{i}\cap C_{j}$ in
$\sigma $.  

Put $b_{i,j} = |A_{i}\cap B_{j}|$ and $c_{i,j}=|A_{i}\cap C_{j}|$.
Define $V_{i,j}$ as before, and set $W_{i,j} = v(\sigma ,k)+\chi
(k\leq r_{1})$, where $\sigma _{k}$ is the largest element of
$A_{i}\cap C_{j}$.  Then similar reasoning yields the formula
\begin{equation}\label{e:super-fermionic}
\langle D_{n}(z;q,t),e_{\eta } h_{\mu }  \rangle = \sum _{\substack{\sigma \in S_{n}\\
\text{$\tilde{\sigma} $ is a $\mu,\eta $-shuffle}}} t^{\comaj
(\sigma)} (\prod _{i,j} \qbinom{V_{i,j}}{b_{i,j}}_{q}) ( \prod _{i,j}
q^{\binom{c_{i,j}}{2}} \qbinom{W_{i,j}}{C_{i,j}}_{q}) .
\end{equation}

\begin{remarks}
(1) In the above, we do not assume that $\mu $ and $\eta $ are
partitions.  Thus ``the'' fermionic formula is really a separate
formula for each possible ordering of the parts of $\mu $ and $\eta $.
The symmetry part of Theorem~\ref{thm:symmetry-positivity} is
equivalent to the statement that all these formulas yield the same
result.

(2) The reader may enjoy verifying that the special case $\mu
=\emptyset $, $\eta =(n)$ of \eqref{e:super-fermionic} agrees with the
fermionic formula in \cite{Haglund03a} for the $q,t$-Catalan
polynomial.
\end{remarks}

\subsection{Schr\"oder paths}
\label{s:schroder}

In \cite{EgHaKiKr03}, Egge, Haglund, Killpatrick and Kremer
conjectured a combinatorial formula for $\langle \nabla e_n
,e_{n-d}h_d \rangle$.  In this section we first show how this
conjecture is a special case of \eqref{e:super-main-equiv}.  Then we
briefly discuss how some ideas in Haglund's recent proof of their
conjecture suggest a refinement of Conjecture \ref{conj:main}.

A {\it Schr\"oder path} is a lattice path from $(n,0)$ to $(0,n)$
composed of steps of the form $(-1,0)$ (south), $(0,1)$ (east) and
$(-1,1)$ (diagonal) which never goes above the line
$i+j=n$.\footnote{In \cite{EgHaKiKr03}, partitions are drawn in the
fourth quadrant, English style, where we draw them in the first
quadrant, French style.  Correspondingly, our Schr\"oder paths are
mirror images of those in \cite{EgHaKiKr03}.}  Egge, {\it et al.}\
gave two different formulations of their conjecture, one involving a
pair of statistics $(\area ,\bounce )$ on Schr\"oder paths and another
involving a pair of statistics $(\dinv ,\area )$.  They showed that
the two formulations are equivalent by exhibiting a bijection which
sends $(\dinv ,\area )$ to $(\area ,\bounce )$.

Given a Schr\"oder path $\Pi $, let $\lambda (\Pi)$ be the partition
whose associated Dyck path is obtained by replacing each diagonal step
of $\Pi$ by a south step followed by an east step.  If $\Pi$ has $d$
diagonal steps, let $T(\Pi )$ be the super tableau of shape $(\lambda
+(1^{n}))/\lambda $ obtained by placing the number $1$ in each square
bordered by one of the $d$ new pairs of south, east steps replacing
the former diagonal steps, and setting all other entries to $\bar{1}$.
The reader will have no problem checking from the definitions of
$\dinv (\Pi)$ and $\area (\Pi )$ given in \cite{EgHaKiKr03} that
$\dinv (\Pi) = \dinv (T(\Pi))$ and  $\area (\Pi ) = |\delta
_{n}/\lambda (\Pi )|$.  Thus the conjecture of Egge, {\it et al.}\
is equivalent to the special case of \eqref{e:super-main-equiv} for
$e_{\eta }h_{\mu } = e_{n-d}h_{d}$.

More recently, Haglund \cite{Haglund03} has proven the conjecture of
Egge, {\it et al.}, and also the $h_{\mu } = h_{n-d}h_{d}$ case of
\eqref{e:main-equiv}.  An important role in both proofs is played by
functions $E_{n,k}$, defined as the coefficients in the Newton
interpolation series expansion
\begin{equation}\label{33}
e_n[Z \frac {1-u}{1-q}] = \sum _{k=1}^n \frac {(u;q)_k}{(q;q)_k}
E_{n,k}(z).
\end{equation}
Note that by setting $u=q$ in \eqref{33} we see $\sum _{k=1}^n E_{n,k}
= e_n$.  The $E_{n,k}$ were first introduced by Garsia and Haglund in
their proof of \eqref{e:qt-Catalan-theorem}.  In \cite{Haglund03} it is
conjectured that
\begin{equation}\label{44}
\langle \Delta _{s_{\beta}} \nabla E_{n,k}, s_{\mu} \rangle \in
{\mathbb N} [q,t]
\end{equation}
for all $\mu, \beta$.  Here $\Delta _f$ is a linear operator
defined on the modified Macdonald basis as
\begin{equation}
\Delta _f {\tilde H}_{\mu} = f[B_{\mu}] {\tilde H}_{\mu},
\end{equation}
where 
\begin{equation}
B_{\mu} = \sum_{(i,j)\in \mu } t^{i}q^{j}.
\end{equation}
In particular, \eqref{44} implies
\begin{equation}\label{46}
\langle \nabla E_{n,k} ,s_{\mu} \rangle \in {\mathbb N}[q,t].
\end{equation}
The conjectured truth of \eqref{44} is largely motivated by
\cite[Cor.~3.5]{Hai02}, which implies
\begin{equation}\label{55}
\langle \Delta _{s_{\beta}} \nabla e_{n}, s_{\mu} \rangle \in {\mathbb
N} [q,t].
\end{equation}

We now introduce a refinement of Conjecture~\ref{conj:main} which
implies \eqref{46}, namely
\begin{conj}\label{conj:nabla-En,k}
For $1\le k \le n$,
\begin{equation}\label{e:nabla-En,k}
\nabla E_{n,k} =
\sum _{\substack{\lambda \subseteq \delta _{n} \\
|\{i:\lambda _{i} = n-i \}| = k}} t^{|\delta
_{n}/\lambda |}D_n^{\lambda}(z;q).
\end{equation}
\end{conj}

\section{Proof of Theorem \ref{thm:symmetry-positivity}}
\label{llt}

In this section we will summarize some results of Lascoux, Leclerc and
Thibon \cite{LaLeTh97,LecThi00}, adding to these a new description of
their ``spin'' statistic on ribbon tableaux, and use all this to prove
Theorem~\ref{thm:symmetry-positivity}.

\subsection{Cores and quotients}
\label{s:cores-quotients}

We begin by recalling some standard facts from the combinatorial
theory of $n$-cores and $n$-quotients, as developed for instance in
\cite{JamKer81, StWh85}.  An {\it $n$-ribbon} is a connected skew
shape of size $n$ and depth $1$, {\it i.e.}, containing no $2\times 2$
rectangle.  A partition $\mu $ is an {\it $n$-core} if there is no
$\nu \subseteq \mu $ such that $\mu /\nu $ is an $n$-ribbon.  Every
$\mu $ contains a unique $n$-core $\nu = \core _{n}(\mu )$ such that
$\mu /\nu $ can be tiled by $n$-ribbons.  In other words, if we
successively remove as many $n$-ribbons from $\mu $ as possible, the
shape $\nu $ that remains does not depend on any choices made.

The {\it content} of a cell $x = (i,j)\in \NN \times \NN $ is defined
as $c(x) = j-i$.  Define the content of an $n$-ribbon to be the
maximum of the contents of its cells.  If $\nu $ is an $n$-core, there
are exactly $n$ shapes $\mu $ such that $\mu /\nu $ is an $n$-ribbon,
and the contents $s_{0},s_{1},\ldots,s_{n-1}$ of these ribbons are
distinct $(\mmod\, n)$.  We always index them so that
$s_{i}\equiv i \pmod{n}$.  Then the $n$-cores are in one-to-one
correspondence with cosets $(s_{0},s_{1},\ldots,s_{n-1})+\ZZ \cdot
(n,n,\ldots,n)$ satisfying $s_{i}\equiv i \pmod{n}$ for all $i$.

Fix an $n$-core $\nu $ with content sequence
$(s_{0},s_{1},\ldots,s_{n-1})$.  If $(\mu ^{(0)},\mu ^{(1)},\ldots,\mu
^{(n-1)})$ is any $n$-tuple of partitions, we define the {\it adjusted
content} of a cell $x\in \mu ^{(i)}$ to be
\begin{equation}\label{e:adj-content}
\tilde{c}(x) = nc(x)+s_{i}.
\end{equation}
Note that $\tilde{c}(x)$ determines which $\mu ^{(i)}$ the cell $x$
belongs to via its congruence class $(\operatorname{mod}\, n)$.  Let
$\Pcal $ be the set of all partitions, and let $\Pcal _{\nu } = \{\mu
\in \Pcal : \core _{n}(\mu ) = \nu \}$.  There is a bijection
\begin{equation}\label{e:quot-n}
\quot _{n}\colon \Pcal _{\nu }\rightarrow \Pcal ^{n},
\end{equation}
written $\quot _{n}(\mu ) = (\mu ^{(0)},\mu ^{(1)},\ldots,\mu
^{(n-1)})$, characterized by the following property: in any tiling of
$\mu /\nu $ by $n$-ribbons, the multiset of contents of the ribbons is
equal to the multiset of adjusted contents $\tilde{c}(x)$, taken over
all $i$ and all cells $x\in \mu ^{(i)}$.  In particular, we have
\begin{equation}\label{e:quot-size}
|\mu /\nu | = n\left|\quot _{n}(\mu )\right| \defeq n\sum _{i}|\mu ^{(i)}|.
\end{equation}
For $\lambda ,\mu \in \Pcal _{\nu }$ we also have 
\begin{equation}\label{e:quot-preserves-containment}
\lambda \subseteq \mu \quad \Leftrightarrow \quad \text{$\lambda
^{(i)}\subseteq \mu ^{(i)}$ for all $i$}.
\end{equation}
Therefore, $\quot _{n}$ extends to a bijection $\quot _{n}(\mu
/\lambda ) = ( \mu ^{(0)}/\lambda ^{(0)},\ldots, \mu ^{(n-1)}/\lambda
^{(n-1)})$ from skew shapes $\mu /\lambda $ with $\lambda ,\mu \in
\Pcal _{\nu }$ to $n$-tuples of skew shapes.  To avoid notational
ambiguity, we henceforth apply $\quot _{n}$ only to skew shapes, so if
$\mu $ is a partition with $\core _{n}(\mu ) = \nu $, we would write $
(\mu ^{(0)},\mu ^{(1)},\ldots,\mu ^{(n-1)}) = \quot _{n}(\mu /\nu )$.

We remark that a given skew shape may have multiple representations
$\theta = \mu /\lambda $ with different $n$-cores $\nu =\core _{n}(\mu
) = \core _{n}(\lambda )$.  However, the resulting $n$-quotients
$\quot _{n}(\theta )$ differ only by translations of the components $\theta
^{(i)}$, which compensate for the change in the contents $s_{i}$
associated with $\nu $ so that the adjusted contents $\tilde{c}(x)$
for $x\in \quot _{n}(\theta )$ remain the same.

A standard, semistandard or super tableau on an $n$-tuple of shapes
$(\mu ^{(0)},\ldots,\mu ^{(n-1)})$ just means a tableau of the
specified sort on the disjoint union of the shapes $\mu ^{(i)}$.  A
{\it standard $n$-ribbon tableau} on a (skew) shape $\mu $ is a tiling
of $\mu $ by $n$-ribbons and a function $T\colon \mu \rightarrow
\{1,2,\ldots,N = |\mu |/n \}$, weakly increasing on each row and
column, which is constant on each ribbon and induces a bijection from
the ribbons to $\{1,2,\ldots,N \}$.  It follows from
\eqref{e:quot-preserves-containment} that $\quot _{n}$ induces a
bijection between standard $n$-ribbon tableaux of shape $\mu $ and
$\SYT (\quot _{n}(\mu ))$.

Call $\mu $ a {\it horizontal} (resp.\ {\it vertical}) {\it $n$-ribbon
strip} if it can be tiled by $n$-ribbons and each $\mu ^{(i)}$ is a
horizontal (resp.\ vertical) strip.  Then there is a unique standard
ribbon tableau $T$ of shape $\mu $ in which the ribbons are labelled
in increasing (resp.\ decreasing) order of content---namely, $\quot
_{n}(T)$ is the unique standard tableau of shape $\quot _{n}(\mu )$ in
which the cells are labelled in increasing (resp.\ decreasing) order
of adjusted content $\tilde{c}(x)$.  The ribbon tiling given by this
distinguished tableau is the {\it official tiling} of the strip $\mu
$.  In the case of a horizontal ribbon strip $\mu $, the official
tiling is characterized by the property that the cell of maximum
content in each ribbon is the minimum cell in its column in the shape
$\mu $.

A {\it semistandard $n$-ribbon tableau} of shape $\mu $ is a tiling of
$\mu $ by $n$-ribbons and a function $T\colon \mu \rightarrow \Acal
_{+}$, weakly increasing on each row and column and constant on each
ribbon, such that for each $a\in \Acal _{+}$, $T^{-1}(a)$ is a
horizontal $n$-ribbon strip with the official tiling.  Define $z^{T}$
to be the product of $z_{T(\theta )}$ over the $n$-ribbons $\theta $
in the given tiling of $\mu $, or equivalently, $\prod _{x\in \mu
}z_{T(x)} = (z^{T})^{n}$.  It is immediate that $\quot _{n}$ induces a
weight-preserving bijection between the set $\SSRT ^{n}(\mu )$ of
semistandard $n$-ribbon tableaux of shape $\mu $ and $\SSYT (\quot
_{n}(\mu ))$.

\subsection{Spin generating functions}
\label{s:spin}

As in \cite{LaLeTh97}, the {\it spin} $s(\theta )$ of an $n$-ribbon
$\theta $ is one less than the number of its rows.  Given a
(semi)standard $n$-ribbon tableau $T$, we set $s(T)$ equal to the sum
of $s(\theta )$ over all ribbons $\theta $ in the tiling underlying
$T$.  One proves that $(-1)^{s(T)}$ depends only on the shape of $T$.
Let $\smin (\mu )$ and $\smax (\mu )$ be the minimum and maximum of
$s(T)$ over all $T$ of shape $\mu $.  The {\it spin} and {\it cospin}
of a tableau $T$ of shape $\mu $ are then defined to be the integers
$\spin (T) = \frac{1}{2}(s(T)-\smin (\mu ))$ and $\cospin (T) =
\frac{1}{2}(\smax (\mu )-s(T))$, respectively.

\begin{thm}[\cite{LaLeTh97}]\label{thm:llt-sym}
For every (skew) shape $\mu $, the spin generating function
\begin{equation}\label{e:G-mu}
G_{\mu }(z;q) = \sum _{T\in \SSRT ^{n}(\mu )} q^{\spin (T)}z^{T}
\end{equation}
is a symmetric function.
\end{thm}

To apply this theorem in our setting, we need an alternative
description of spin.  Fix a content sequence
$(s_{0},s_{1},\ldots,s_{n-1})$ with $s_{i} \equiv i\pmod{n}$ and an
$n$-tuple of shapes $\boldmu = (\mu ^{(0)},\mu ^{(1)},\ldots,\mu
^{(n-1)})$.  Let $S$ be a semistandard tableau of shape $\boldmu $.
An {\it inversion} is a pair of entries $S(x)=a$, $S(y)=b$ such that
$a<b$ and $0<\tilde{c}(x)-\tilde{c}(y)<n$.  Denote by $\inv (S)$ the
number of inversions in $S$.

\begin{lemma}\label{lemma:new-spin-special}
Given a (skew) shape $\mu $, there is a constant $e$ such that for
every standard $n$-ribbon tableau $T$ of shape $\mu $, we have
$\spin (T) = e-\inv (\quot _{n}(T))$.
\end{lemma}

\begin{figure}
\begin{center}
\includegraphics{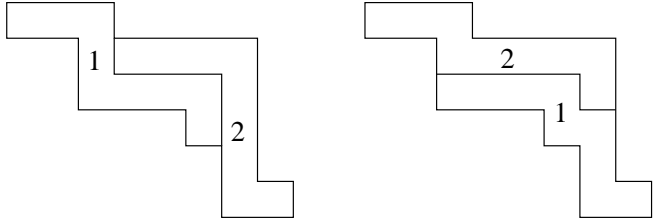}
\caption{\label{fig:two-ribbons}%
A shape with two tilings by two $n$-ribbons}
\end{center}
\end{figure}

\begin{proof}
Say that standard tableaux $S$, $S'$ of shape $\boldmu =\quot _{n}(\mu
)$ differ by a {\it switch} if they are identical except for the
positions of two consecutive entries $a$, $a+1$.  Every pair of
tableaux is connected by a sequence of switches.  If $S$ and $S'$
differ by a switch, then the only entries that might form an inversion
in one tableau but not the other are $a$, $a+1$.  Assuming as we may
that $\inv (S')\geq \inv (S)$, we therefore have $\inv (S') = \inv
(S)+\epsilon $ with $\epsilon \in \{0,1 \}$.

Let $S = \quot _{n}(T)$, $S' = \quot _{n}(T')$.  What must be shown is
that $\spin (T) = \spin (T')+\epsilon $.  Since $T$ and $T'$ are
identical except for the ribbons labelled $a$ and $a+1$, the problem
reduces to the case $|\mu |=2n$, $|\boldmu | =2$.  Let $x$, $y$ be the
two cells of $\boldmu$, with $S(x)=1$, $S(y)=2$ and $S'(x)=2$,
$S'(y)=1$.  Note that $\tilde{c}(x) \not = \tilde{c}(y)$.  We have
$\epsilon =0$ if and only if $|\tilde{c}(y)-\tilde{c}(x)|\geq n$.
This means that $\mu $ is the union of two ribbons whose cells have no
contents in common, and these ribbons are unique.  Hence $\spin (T) =
\spin (T')$.

Conversely, we have $\epsilon =1$ if and only if
$|\tilde{c}(y)-\tilde{c}(x)|< n$, which means that $\mu $ can be tiled
by two ribbons whose cells have at least one content in common.  In
this case $\mu $ has exactly two ribbon tilings, as in
Figure~\ref{fig:two-ribbons}, each supporting one of the standard
tableaux $T$, $T'$.  In fact, $\mu $ is both a horizontal and a
vertical $n$-ribbon strip, the tiling of $T$ is the official tiling of
$\mu $ as a horizontal strip, and that of $T'$ is the official tiling
as a vertical strip.  In the horizontal tiling, each ribbon has one
more row than the corresponding ribbon with the same content in the
vertical tiling (see \cite[Lemma~4.1]{Pak00}).  Hence $\spin (T)=\spin
(T')+1$.
\end{proof}

\begin{lemma}\label{lemma:new-spin-general}
Lemma~\ref{lemma:new-spin-special} holds also for semistandard tableaux.
\end{lemma}

\begin{proof}
If $S$ is a standard tableau of shape $\boldmu = \quot _{n}(\mu )$,
call $a$ a {\it descent} of $S$ if $S(x) = a$, $S(y)=a+1$, with
$\tilde{c}(x)>\tilde{c}(y)$.  If $\boldmu $ is a horizontal strip,
there is a unique standard tableau of shape $\boldmu $ with no
descents, and conversely, such a tableau exists only on a horizontal
strip.  As in the proof of Proposition~\ref{prop:super-schur-via-Q},
it follows that each semistandard tableau $T$ has a unique standardization
$S$ such that $T\circ S^{-1}$ is weakly increasing, and if $T\circ
S^{-1}(j)=T\circ S^{-1}(j+1)=\cdots =T\circ S^{-1}(k)$, then $d(T)\cap
\{j,\ldots,k-1 \} = \emptyset $.

By definition, equal entries $T(x)=T(y)=a$ contribute nothing to $\inv
(T)$.  In the standardization $S$, equal entries are replaced with
entries labelled in increasing order of adjusted content $\tilde{c}$,
which contribute nothing to $\inv (S)$.  On the other hand, unequal
entries of $T$ give rise to entries ordered in the same way in $S$, so
$\inv (T) = \inv (S)$.

Given a semistandard $n$-ribbon tableau $T$, we define its
standardization to be the unique standard $n$-ribbon tableau $S$ such
that $\quot _{n}(S)$ is the standardization of $\quot _{n}(T)$, as
above.  The important point to notice here is that for each letter
$a$, the horizontal ribbon strip $T^{-1}(\{a \})$ has its ribbons
labelled in $S$ in increasing order of content, and hence it is tiled
in $S$ by the official tiling.  This shows that $T$ and its
standardization have the same underlying ribbon tiling, and hence the
same spin $\spin (T) = \spin (S)$.

These observations reduce the lemma for semistandard tableaux to the
standard case.
\end{proof}

\begin{remark}\label{rem:Schilling}
Schilling, Shimozono and White \cite{SchShiWh03} defined an inversion
number $\inv '(T)$ such that $\cospin (T) = \inv '(\quot (T))$
exactly, without the constant error term $e$ in
Lemma~\ref{lemma:new-spin-special}.  An inversion by their definition
is an inversion by ours which also satisfies some extra conditions.
Haiman's student Michelle Bylund and Haiman found the simpler
definition used here and the proof given above.
\end{remark}

\begin{cor}\label{cor:Gmu-by-inv}
Fix an $n$-tuple of shapes $\boldmu =(\mu ^{(0)},\ldots,\mu ^{(n-1)})$
and a sequence of content offsets $s_{i}\equiv i\pmod{n}$ as in the
definitions of $\tilde{c}$ and $\inv $.  Then there exists $\mu $ with
$\quot _{n}(\mu ) = \boldmu $ such that
\begin{equation}\label{e:Gmu-by-inv}
q^{e}G_{\mu }(z;q^{-1}) = \sum _{T\in \SSYT (\boldmu )} q^{\inv
(T)} z^{T}
\end{equation}
for some exponent $e$.  In particular, the expression on the
right-hand side is a symmetric function.
\end{cor}

\begin{proof}
This is immediate from Lemma~\ref{lemma:new-spin-general} if
$(s_{0},s_{1},\ldots,s_{n-1})$ is the content sequence of some
$n$-core $\nu $.  But we can always make it one by adding
$c(n,n,\ldots,n)$ for some $c$.  This change does not alter the value
of $\inv (T)$.
\end{proof}

\begin{proof}[Proof of Theorem~\ref{thm:symmetry-positivity}
(symmetry)]
We are to show that the expression $D_{n}^{\lambda }(z;q)$ in
\eqref{e:lambda-term} is a symmetric function.  Write $\lambda
=(0^{\alpha _{0}},1^{\alpha _{1}},\ldots , (n-1)^{\alpha _{n-1}})$,
defining $\alpha _{0}$ so that $\sum _{j}\alpha _{j} = n$.  Set $\mu
^{(j)} = (1^{\alpha _{j}})$, a single column of the same height as
column $j$ in $(\lambda +(1^{n}))/\lambda $.  Translating the columns
of $(\lambda +(1^{n}))/\lambda $ onto the corresponding columns $\mu
^{(j)}$ gives a bijection between the cells of $(\lambda
+(1^{n}))/\lambda $ and those of $\boldmu =(\mu ^{(0)},\ldots,\mu
^{(n-1)})$.  This induces in the obvious way a bijection between
semistandard tableaux of these shapes.

Take the offsets $s_{j}= j-n(j+\lambda '_{j+1})$.  Then if $x\in (\lambda
+(1^{n}))/\lambda $ corresponds to $x'\in \mu ^{(j)}$, we have 
\begin{equation}\label{e:c-vs-d}
\tilde{c}(x') = j-nd(x).
\end{equation}
It follows that conditions (i) or (ii) in the definition of
d-inversion for cells $x,y\in (\lambda +(1^{n}))/\lambda $ are
equivalent to the condition $0<\tilde{c}(x')-\tilde{c}(y')<n$ for the
corresponding cells $x',y'\in \boldmu $.  If $T\in \SSYT (\lambda
+(1^{n})/\lambda )$ corresponds to $T'\in \SSYT (\boldmu )$ we
therefore have $\dinv (T) = \inv (T')$.  Hence $D_{n}^{\lambda
}(z;q)$ coincides with the right-hand side of \eqref{e:Gmu-by-inv}
for this choice of $\boldmu $ and $s_{i}$.
\end{proof}

\begin{example}\label{ex:symmetry-proof}
We make the constructions in the proof more explicit for the tableau
$T$ shown in Figure~\ref{fig:park}.  For this $T$, we have $\boldmu =
((1^{2}),\,(1^{2}),\,\emptyset ,\,(1),\,\emptyset
,\,(1),\,(1^{2}),\,\emptyset )$ and
\[
T' = \bigl(\; \tableau{8\\6}\, ,\; \tableau{7\\5}\, ,\; \emptyset \, ,\;
\tableau{1}\, ,\; \emptyset \, ,\; \tableau{2}\, ,\; \tableau{4\\3}\, ,\;
\emptyset\, \bigr).
\]
The content offsets are $(-48, -39, -46, -45, -52, -51, -42, -49)$,
giving adjusted contents
\[
\bigl(\; \tableau{{\scriptstyle -56}\\{\scriptstyle -48}}\, ,\;
\tableau{{\scriptstyle -47}\\{\scriptstyle -39}}\, ,\; \emptyset \,
,\; \tableau{{\scriptstyle -45}}\, ,\; \emptyset \, ,\;
\tableau{{\scriptstyle -51}}\, ,\; \tableau{{\scriptstyle
-50}\\{\scriptstyle -42}}\, ,\; \emptyset\, \bigr).
\]
The reader can verify that each $d$-inversion in $T$ listed
after eq.~\eqref{e:dinv} corresponds to a pair of entries
$T'(x')<T'(y')$ in cells $x'$, $y'$ with adjusted contents
$0<\tilde{c}(x')-\tilde{c}(y')<8$.  For example, the $d$-inversion
$(4,1)$ has $\tilde{c}(x') = -45$, $\tilde{c}(y') = -50$.
\end{example}

\subsection{Positivity}
\label{s:positivity}

When $\mu $ is a partition with $\core _{n}(\mu ) = \emptyset $,
Leclerc and Thibon \cite{LecThi00} have shown that the coefficient
\begin{equation}\label{e:s-lambda-Gmu}
\langle q^{\smin (\mu )} G_{\mu }(z;q^{2}),s_{\lambda }(z) \rangle
\end{equation}
(which is a $q$-analog of the Littlewood-Richardson coefficient
$c^{\lambda }_{\mu ^{(0)},\ldots,\mu ^{(n-1)}}$) coincides with a
parabolic Kazhdan-Lusztig polynomial $P^{-}_{\mu +\rho , n\lambda
+\rho }(q)$ for a suitable affine symmetric group $\widehat{S}_{r}$.
In turn, Kashiwara and Tanisaki \cite{KashTan02} have interpreted the
coefficients of these polynomials as decomposition multiplicities for
certain non-irreducible perverse sheaves on affine partial flag
varieties, which shows that they are positive.  This together with our
work in \S \ref{s:spin} would immediately imply the positivity
part of Theorem~\ref{thm:symmetry-positivity}, were it not for the
fact that we had to introduce non-trivial offsets $s_{i}$.  To account
for them and complete the proof of
Theorem~\ref{thm:symmetry-positivity} we need the following small
improvement on the results of \cite{LecThi00}.

\begin{prop}\label{prop:KL-poly}
Let $\mu $ be a partition and set $\nu = \core _{n}(\mu )$.  Then
\begin{equation}\label{e:Gmu=KL-poly}
\langle q^{\smin (\mu /\nu )} G_{\mu /\nu }(z;q^{2}),s_{\lambda }(z) \rangle =
P^{-}_{\mu +\rho ,\nu +n\lambda +\rho }(q)
\end{equation}
is a parabolic Kazhdan-Lusztig polynomial---written here using the
notation of \cite{LecThi00}.
\end{prop}

Since this is essentially a result of Leclerc and Thibon, we will
confine ourselves to brief remarks on what is needed to deduce it from
the contents of \cite{LecThi00}.  Leclerc and Thibon work in a
$q$-Fock space $\Fcal _{r}$, where $r$ is arbitrary, provided it is
greater than or equal to the length of all partitions under
discussion.  The space $\Fcal _{r}$ is equipped with a natural basis,
whose elements are denoted $\ket{\mu +\rho }$ (with $\rho =\delta
_{r}$), and a Kazhdan-Lusztig type basis, denoted $G^{-}_{\mu +\rho
}$.  The Kazhdan-Lusztig polynomials $P^{-}_{\lambda +\rho ,\mu +\rho
}(-q^{-1})$ are the coefficients of $G^{-}_{\mu +\rho }$ with respect
to the natural basis.  The use of $-q^{-1}$ is an artifact of the
notation; it gets replaced by $q$ in the end.  The algebra of
symmetric functions in $r$ variables acts on $\Fcal _{r}$ in such a
way that $\langle (-q)^{-\smin (\mu/\nu )} G_{\mu /\nu
}(z;q^{-2}),s_{\lambda }(z) \rangle$ is the coefficient of $\ket{\mu
+\rho}$ in $s_{\lambda }\cdot \ket{\nu +\rho }$.

A partition $\nu $ is called {\it $n$-restricted} if $\nu _{i}-\nu
_{i+1}< n$ for all $i$.  Leclerc and Thibon prove that if $\nu $ is
$n$-restricted, then $s_{\lambda }\cdot G^{-}_{\nu+\rho } = G^{-}_{\nu
+n\lambda +\rho }$.  Proposition~\ref{prop:KL-poly} follows
immediately from the foregoing and the following two additional facts.

\begin{lemma}\label{lemma:core-restricted}
Any $n$-core is $n$-restricted.
\end{lemma}

\begin{proof}
Obvious.
\end{proof}

\begin{lemma}\label{lemma:core-is-KL}
If $\nu $ is an $n$-core, then $G^{-}_{\nu +\rho } = \ket{\nu +\rho
}$.
\end{lemma}

\begin{proof}
Using the notation of \cite{LecThi00}, it suffices to prove that
$\overline{\ket{\nu +\rho }} = \ket{\nu +\rho }$.  By
\cite[Prop~5.9]{LecThi00}, there are well-defined elements
$\ket{\gamma }\in \Fcal _{r}$ for any sequence $\gamma =(\gamma
_{1},\ldots,\gamma _{r})$, not necessarily decreasing, satisfying the
following straightening relations when  $\gamma
_{i+1}-\gamma _{i}= kn+j$ with $k\geq 0$ and $0\leq j<n$.
\begin{align}
\ket{\gamma } &= 0\quad \text{if $j=k=0$,}\\
\ket{\gamma } &= -\ket{\sigma_{i}\gamma } \quad \text{if $j=0$ and
$k\not =0$,}\\ 
\ket{\gamma } &= -q^{-1}\ket{\sigma_{i}\gamma } \quad \text{if $k=0$
and $j\not =0$,}\\
\ket{\gamma } &= -q^{-1}\ket{\sigma_{i}\gamma} -
\ket{y_{i}^{-k}y_{i+1}^{k}\gamma} -
q^{-1}\ket{y_{i}^{k}y_{i+1}^{-k}\sigma_{i}\gamma} \quad
\text{otherwise.}
\end{align}
Here $\sigma_{i}$ is the transposition exchanging $\gamma _{i}$ and
$\gamma _{i+1}$ and $y_{i}$ is the operator that adds $n$ to $\gamma
_{i}$.  Let $w_{0}$ denote the longest permutation in $S_{r}$.  By
\cite[Prop.~5.7 and Cor.~5.10]{LecThi00}, we have
\begin{equation}\label{e:ket-bar-expand}
\overline{\ket{\nu
+\rho }} = (-1)^{l(w_{0})}q^{e}\ket{w_{0}(\nu +\rho )} = \ket{\nu +\rho
}+\sum _{\lambda }a_{\lambda +\rho ,\nu +\rho }(q)\ket{\lambda +\rho
}
\end{equation}
for a suitable exponent $e$, and all terms $\ket{\lambda +\rho }$ in
the sum on the right, which arise from the process of straightening
$\ket{w_{0}(\nu +\rho )}$, are lexicographically less than $\nu +\rho
$.  From the straightening relations we see that all these terms have
$|\lambda |=|\nu |$, and the multiset of congruence classes $\lambda
_{i}+\rho _{i}\pmod{n}$ is the same as that of $\nu +\rho $.  However,
this last condition implies that $\core _{n}(\lambda ) = \core
_{n}(\nu )=\nu $, which is absurd, since $|\lambda |=|\nu |$ and
$\lambda \not =\nu $.  In short, the final sum on the right in
\eqref{e:ket-bar-expand} vanishes, yielding the desired result.
\end{proof}

\begin{remark}
It is conjectured that $G_{\mu /\nu }(z;q)$ is Schur positive even
when $\nu $ is not an $n$-core.  One can use the results of
\cite{LecThi00} to write $\langle q^{\smin (\mu /\nu )} G_{\mu /\nu
}(z;q^{2}),s_{\lambda }(z) \rangle$ explicitly in terms of
Kazhdan-Lusztig polynomials, but in general the resulting expressions
contain negative terms.
\end{remark}

\begin{proof}[Proof of Theorem~\ref{thm:symmetry-positivity}
(positivity)]
In the proof of the symmetry part of the theorem, we
have shown that $D_{n}^{\lambda }(z;q) = q^{e}G_{\mu }(z;q^{-1})$ for
a suitable exponent $e$ and skew shape $\mu $.  For this $\mu $, each
$\mu ^{(i)} = (1^{\alpha _{i}})$ is a partition shape, so $\mu = \eta
/\nu$, where $\nu = \core _{n}(\eta )$ is the $n$-core associated to
the specified offsets $s_{i}$.  Hence $G_{\mu }(z;q)$ is Schur
positive by Proposition~\ref{prop:KL-poly}.
\end{proof}

\section{Higher powers}
\label{higher}

In this section we explore to what extent the preceding conjectures
and results generalize to the higher powers $\nabla ^{m}e_{n}(z)$.

\subsection{The meaning of $\nabla ^{m}e_{n}(z)$}
\label{s:meaning}

As explained in the introduction, $\nabla e_{n}(z)$ is the Frobenius
series of the diagonal coinvariant ring $R_{n}$.  The higher powers
$\nabla ^{m}e_{n}(z)$ have a similar interpretation, which shows that
they are also Schur positive.

\begin{prop}\label{prop:nabla-m=Frobenius}
Let $I = ((\xx ,\yy )\cap \CC [\xx ,\yy ]^{S_{n}})$ be the ideal in
$\CC [\xx ,\yy ]$ generated by all $S_{n}$-invariant polynomials
without constant term, and  $J = (\CC [\xx ,\yy ]^{\epsilon })$ the
ideal generated by all antisymmetric polynomials.  Then
\begin{equation}\label{e:nabla-m=Frobenius}
\nabla ^{m}e_{n}(z) = \Fcal _{\varepsilon ^{m-1} \otimes
J^{m-1}/IJ^{m-1}}(z;q,t)
\end{equation}
where $\varepsilon $ is the sign representation.
\end{prop}

\begin{proof}
Although this follows from the methods of \cite{Hai02}, it wasn't
shown explicitly there, so we explain what more needs to be said.  As
in \cite[Cor.~3.5, eq.~(110)]{Hai02}, the quantity $\nabla
^{m}e_{n}(z)$ coincides with \cite[Theorem~3.3, eq.~(107)]{Hai02},
with the factor $s_{\nu }[B_{\mu }(q,t)]$ there replaced by
$e_{n}^{m-1}[B_{\mu }(q,t)] = t^{(m-1)n(\mu )}q^{(m-1)n(\mu ')}$.  It
follows from that theorem that $\nabla ^{m}e_{n}(z)$ is the Frobenius
series of
\begin{equation}\label{e:complicated-polygraph}
(R(n,(m-1)n)/{\mathfrak m}R(n,(m-1)n))^{\epsilon },
\end{equation}
where $l=(m-1)n$, $R(n,l)$ and ${\mathfrak m}$ are as in \cite{Hai02},
and $(-)^{\epsilon }$ denotes the space of antisymmetric elements with
respect to the action of $S_{n}^{m-1}\subseteq S_{l}$.  The proof of
\cite[Prop.~4.11.1]{Hai01} identifies $R(n,(m-1)n)^{\epsilon }$ with
$\varepsilon ^{m-1}\otimes J^{m-1}$, and with this identification,
\eqref{e:complicated-polygraph} becomes $\varepsilon ^{m-1}\otimes 
J^{m-1}/IJ^{m-1}$ in our present notation.

\end{proof}

\subsection{An extension of Conjecture~\ref{conj:main}}
\label{s:extended-conjecture}

We begin by generalizing the notion of d-inversion and the statistic
$\dinv (T)$.  Our new definitions depend on the integer $m$ and reduce
for $m=1$ to the definitions of d-inversion and $\dinv (T)$ in \S
\ref{s:main}.

For each cell $x=(i,j)\in \NN \times \NN $, put $d_m(x) = mi+j$ (this
keeps track of which diagonal of slope $-1/m$ contains $x$).  Given
$\lambda \subseteq m \delta _{n}$, let $T$ be a semistandard tableau
of shape $(\lambda +(1^{n}))/\lambda $.  Let $T(x) = a$, $T(y)= b$ be
two entries with $a<b$, and put $x=(i,j)$, $y=(i',j')$.  We say that
this pair of entries contributes anywhere from $0$ to $m$ {\it
d-inversions}, according to the following rules:
\begin{itemize}
\item [(i)] if $j>j'$, this pair contributes $\max (0,m-|d_m(y)-d_m(x)|)$
inversions; 
\item [(ii)] if $j<j'$, it contributes $\max (0,m-|d_m(y)-d_m(x)-1|)$
inversions. 
\end{itemize}
We also allow equal entries to contribute.  Define the reverse
diagonal lexicographic order by
\[
x <_{d} y \quad \text{if $d_m(x)>d_m(y)$, or $d_m(x)=d_m(y)$ and $j<j'$,
where $x=(i,j)$, $y=(i',j')$}.
\]
Then a pair of equal entries $T(x)=T(y)=a$ contributes the same number
of d-inversions as would a pair of unequal entries $T(x)=a$, $T(y)=b$
with $x <_{d} y$ and $a<b$.

We extend these definitions to super tableaux by applying the rules
above for unequal entries $T(x)=a$, $T(y)=b$ and for equal entries
$T(x)=T(y)=a\in \Acal _{+}$.  For negative entries, we use the
opposite rule: a pair of entries $T(x)=T(y)=a\in \Acal _{-}$
contributes the same number of d-inversions as would a pair of unequal
entries $T(x)=a$, $T(y)=b$ with $x >_{d} y$ and $a<b$.

\begin{remark}
Another way to formulate the rule for a pair of equal entries
$T(x)=T(y)=a$ is that if $a$ is positive, it contributes the minimum
of the two alternatives described by (i) and (ii) above; if $a$ is
negative, it contributes the maximum.
\end{remark}

Let $\dinv _{m}(T)$ denote the total number of d-inversions
contributed by pairs of entries in $T$.  The extensions of
Definition~\ref{defn:Dn} and Conjecture~\ref{conj:main} are as
follows.

\begin{defn}\label{defn:Dnm}
\begin{equation}\label{e:Dnm}
D_{n}^{(m)}(z;q,t) = \sum _{\lambda \subseteq m\delta _{n}}\; \sum
_{T\in \SSYT (\lambda +(1^{n})/\lambda )} t^{|m\delta _{n}/\lambda |}
q^{\dinv _{m} (T)} z^{T}.
\end{equation}
\end{defn}

\begin{conj}\label{conj:m}
We have the identity
\begin{equation}\label{e:m}
\nabla ^{m} e_{n}(z) = D_{n}^{(m)}(z;q,t).
\end{equation}
Equivalently, for all $\mu $ we have 
\begin{equation}\label{e:m-equiv}
\langle \nabla ^{m} e_{n}, h_{\mu } \rangle = \sum _{\lambda \subseteq
m \delta _{n}}\; \sum _{T\in \SSYT (\lambda +(1^{n})/\lambda ,\, \mu )}
t^{|m \delta_{n} / \lambda|} q^{\dinv _{m} (T)}.
\end{equation}
\end{conj}

Theorem~\ref{thm:symmetry-positivity} now generalizes just as we should expect.

\begin{thm}\label{thm:m-symmetry}
The quantity $D_{n}^{(m)}(z;q,t)$ is a symmetric function in $z$, and
it is Schur positive.  In fact, each term
\begin{equation}\label{e:m-lambda-term}
D_{n}^{(m),\lambda }(z;q) = \sum _{T\in \SSYT (\lambda
+(1^{n})/\lambda )} q^{\dinv _{m} (T)} z^{T}.
\end{equation}
is individually symmetric and Schur positive.
\end{thm}

Theorem~\ref{thm:m-symmetry} will be proven in \S \ref{s:m-symmetry}.
Granting it for the moment, let us deduce some consequences.

\begin{thm}\label{thm:super-Dnm}
The superization $\tilde{D}_{n}^{(m)}(z,w;q,t) = \omega ^{W}
D_{n}^{(m)}[Z+W; q,t]$ is given by
\begin{equation}\label{e:super-Dnm}
\tilde{D}_{n}^{(m)}(z,w;q,t) = \sum _{\lambda \subseteq m\delta
_{n}}\; \sum _{T\in \SSYT _{\pm }(\lambda +(1^{n})/\lambda )}
t^{|m\delta _{n}/\lambda |} q^{\dinv _{m} (T)} z^{T}.
\end{equation}
Equivalently, Conjecture~\ref{conj:m} implies
\begin{equation}\label{e:super-Dnm-equiv}
\langle \nabla ^{m} e_{n}, e_{\eta }h_{\mu } \rangle = \sum _{\lambda
\subseteq m \delta _{n}}\; \sum _{T\in \SSYT_{\pm } (\lambda
+(1^{n})/\lambda ,\, \mu, \,\eta )} t^{|\delta_{n} / \lambda|}
q^{\dinv _{m} (T)}.
\end{equation}
\end{thm}

\begin{proof}
The proof of Theorem~\ref{thm:super-Dn} applies almost verbatim.  Only
the verification that a super tableau $T$ and its standardization $S$
satisfy $\dinv (S)=\dinv (T)$ needs to be adapted to the case of
general $m$.  But this is immediate, given our rule for the number of
d-inversions contributed by a pair of equal entries.
\end{proof}

\subsection{Specializations}
\label{s:m-specializations}

We examine the analogs for $D_{n}^{(m)}(z;q,t)$ of some of the
specializations of $D_{n}(z;q,t)$ described in \S
\ref{specializations}.  Beginning with the easiest specialization, at
$q=1$, we have by \cite[Thms.~5.2, 5.3]{GaHa96} the following analog
of \eqref{e:nabla-en-q=1}:
\begin{equation}\label{e:nabla-m-en-q=1}
\nabla ^{m} e_{n}(z)|_{q=1} = \sum _{\lambda \subseteq m\delta _{n}}
t^{|m\delta _{n}/\lambda |} e_{\alpha }(z),
\end{equation}
where $\lambda =(0^{\alpha _{0}},1^{\alpha _{1}},2^{\alpha
_{2}},\ldots )$ and $\sum _{i}\alpha _{i} = n$.  Clearly, this
coincides with $D_{n}^{(m)}(z;1,t)$.

For the specializations $t=0$ and $q=0$, we first observe that the
analog of Lemma~\ref{lemma:t=0} is the identity
\begin{equation}\label{e:m-nabla-t=0}
\nabla ^{m}e_{n}(z)|_{t=0} = q^{(m-1)\binom{n}{2}} (q;q)_{n}
h_{n}[Z/(1-q)].
\end{equation}
This can be deduced from Lemma~\ref{lemma:t=0} either by observing that the
right-hand side of \eqref{e:t=0} is the modified Macdonald polynomial
$\tilde{H}_{(n)}$, or by using
Proposition~\ref{prop:nabla-m=Frobenius} and the fact that
$J^{m-1}\cap \CC [\xx ]$ is the principal ideal generated by $\Delta
(\xx )^{m-1}$, where $\Delta (\xx )$ is the Vandermonde determinant in
the variables $\xx =x_{1},\ldots,x_{n}$.  Multiplication by $\Delta
(\xx )^{m-1}$ induces an isomorphism $R_{n}\cap \CC [\xx ]\rightarrow
(J^{m-1}\cap \CC [\xx ])/(IJ^{m-1}\cap \CC [\xx ])$ which is
homogeneous of degree $(m-1)\binom{n}{2}$.

\begin{prop}\label{prop:m-t=0}
We have
\begin{equation}\label{e:Dnm-t=0}
D_{n}^{(m)}(z;q,0) = \nabla^{m} e_{n}(z)|_{t=0}.
\end{equation}
\end{prop}

\begin{proof}
Same as the proof of Proposition~\ref{prop:t=0} except that now the
sole term is $\lambda =m\delta _{n}$, and $\dinv _{m} (T)$ is the number of
ordinary inversions plus $(m-1)\binom{n}{2}$.
\end{proof}

\begin{prop}\label{prop:m-q=0}
We have
\begin{equation}\label{e:Dnm-q=0}
D_{n}^{(m)}(z;0,t) = \nabla^{m} e_{n}(z)|_{q=0}.
\end{equation}
\end{prop}

\begin{proof}
One can verify that the criterion in Lemma~\ref{lemma:dinv=0} for a
standard tableau $T\in \SYT (\lambda +(1^{n})/\lambda )$ to have
$\dinv _{m} (T) = 0$ is the same for general $m$ as it is for $m=1$.  The
only difference for $m>1$ is that $|m\delta _{n}/\lambda | = |\delta
_{n}/\lambda |+(m-1)\binom{n}{2}$, which shows that
$D_{n}^{(m)}(z;0,t) = t^{(m-1)\binom{n}{2}}D_{n}(z;0,t)$.  Exchanging
$q$ and $t$ in \eqref{e:m-nabla-t=0}, and using
Proposition~\ref{prop:q=0}, we see that this agrees with
\eqref{e:Dnm-q=0}.
\end{proof}

Next we turn to the Catalan specialization.  Higher $q,t$-Catalan
polynomials were defined in \cite{GaHa96} (see also \cite{Hai98}) by a
formula which amounts to
\begin{equation}\label{e:m-qt-Catalan}
C_{n}^{(m)}(q,t) = \langle \nabla ^{m}e_{n},e_{n} \rangle.
\end{equation}
From \eqref{e:super-Dnm-equiv}, we see that Conjecture~\ref{conj:m}
implies the following analog of \eqref{e:qt-Catalan-theorem}.
\begin{equation}\label{e:m-qt-Catalan-dinv}
C_{n}^{(m)}(q,t) = \sum _{\lambda \subseteq m \delta _{n}} t^{|m
\delta _{n}/\lambda |}q^{\dinv _{m} (\lambda )}.
\end{equation}
The analog of \eqref{e:qt-Catalan-hook} is a conjecture of Haiman (see
also \cite{Loehr03}) that
\begin{equation}\label{e:m-qt-Catalan-hook}
C_{n}^{(m)}(q,t) = \sum _{\lambda \subseteq m\delta _{n}} t^{|m\delta
_{n}/\lambda | } q^{b_{m}(\lambda )},
\end{equation}
where $b_{m}(\lambda )$ is the number of cells $x\in \lambda $ such that
\begin{equation}\label{e:bm-lambda}
ml(x)\leq a(x)\leq ml(x)+m.
\end{equation}
We will verify that
\eqref{e:m-qt-Catalan-dinv} and \eqref{e:m-qt-Catalan-hook} are
equivalent (but not that they are true---see \S \ref{problems}, Problem~\ref{prob:m-cat}).

\begin{lemma}\label{lemma:bm-lambda}
We have $b_{m}(\lambda ) = \dinv _{m} (T)$, where $T$ is the
super tableau of shape $(\lambda +(1^{n}))/\lambda $ whose every entry
is the negative letter $\bar{1}$.
\end{lemma}

\begin{proof}
Let $\Dcal $ be the lattice path from $(n,0)$ to $(0,mn)$ formed by
the outer boundary of $\lambda $ together with segments of the $i$ and
$j$-axes.  When $m=1$, this is the Dyck path associated with $\lambda
$.  In the general case, it is a lattice path that never goes above
the line $i+j/m = n$.

Given a cell $x\in \lambda $, let $u$ be the cell of $(\lambda
+(1^{n}))/\lambda $ just outside the end of the arm of $x$.  Let $t$ be
the unit-width segment at the end of the leg of $x$.  Project $t$
along diagonals of slope $-1/m$ onto a vertical segment $s$ of $\Dcal
$ of height $1/m$, as indicated in Figure~\ref{fig:uv-projection}, and
let $v$ be the cell of $(\lambda +(1^{n}))/\lambda $ containing $s$ on
its left border.  Note that the point of $\NN \times \NN $ represented
by each cell is its lower-left corner.

\begin{figure}
\begin{center}
\includegraphics{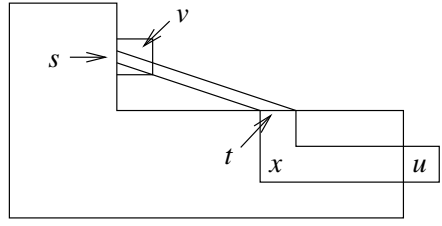}
\caption{\label{fig:uv-projection}%
Construction of the cells $u,v$ corresponding to $x$.}
\end{center}
\end{figure}

If $s_{0}$ is the lower end of the segment $s$ and $t_{0}$ is the left
end of $t$, then we have 
\begin{align}\label{e:d(s0)-vs-d(u)}
d_m(s_{0}) & = d_m(t_{0}) = d_m(x) + ml(x)+m \\
d_m(u) &	= d_m(x) + a(x) +1,
\end{align}
whence $a(x)-ml(x) = d_m(u)-d_m(s_{0})+m-1$.  If $x$ satisfies
\eqref{e:bm-lambda}, we therefore have
\begin{equation}\label{e:d(s0)-range}
d_m(u)-1\leq d_m(s_{0})\leq d_m(u)+m-1.
\end{equation}
Conversely, given $u\in (\lambda + (1^{n}))/\lambda $ and a vertical
segment $s$ of $\Dcal $ of height $1/m$, situated to the left of $u$
and with endpoints in $(\frac{1}{m}\NN)\times \NN $, the inequality
$d_m(s_{0})\leq d_m(u)+m-1$ implies that the diagonal through the upper
endpoint $s_{1}$ lies on or to the left of the upper-left corner of
$u$.  Hence we can project $s$ onto a unit-width horizontal segment
$t$ of $\Dcal $ with $t$ and $u$ belonging to the hook of a (unique)
cell $x\in \lambda $, as in the construction above.  Condition
\eqref{e:bm-lambda} for this $x$ is then equivalent to
\eqref{e:d(s0)-range}.

Now, given $u,v\in (\lambda +(1^{n}))/\lambda $ with $v$ to the left
of $u$, the number of segments $s$ as above that lie on the left
border of $v$ and satisfy \eqref{e:d(s0)-range} is equal to 
\begin{equation}\label{e:number-of-segments}
\min (d_m(v)+m-1,d_m(u)+m-1)-\max (d_m(v),d_m(u)-1)+1,
\end{equation}
or to zero if this expression is negative.  If $d_m(u)>d_m(v)$ this
simplifies to $\max (0,m-(d_m(u)-d_m(v)-1))$.  Otherwise it is $\max
(0,m-(d_m(v)-d_m(u)))$.

In $T$, the cells $u$ and $v$ contain equal
negative entries.  If $d_m(u)>d_m(v)$ then $u<_{d}v$ and this pair
contributes $\max (0,m-|d_m(u)-d_m(v)-1|) = \max (0,m-(d_m(u)-d_m(v)-1))$
d-inversions.  Otherwise, $u>_{d}v$ and it contributes $\max
(0,m-|d_m(v)-d_m(u)|) = \max (0,m-(d_m(v)-d_m(u)))$ d-inversions.  The lemma
is proved.
\end{proof}

The analog for $m>1$ of the Haglund-Loehr conjecture discussed in \S
\ref{s:Haglund-Loehr} is a conjecture of Loehr and Remmel
\cite{LoeRem03} for the value of $\langle \nabla^{m} e_{n},e_{1}^{n}
\rangle$.  It is a simple observation that their conjecture is
equivalent to $\langle \nabla^{m} e_{n},e_{1}^{n} \rangle = \langle
D_{n}^{(m)}(z;q,t),e_{1}^{n} \rangle$.  In this connection we should
mention that Loehr \cite{Loehr03} has given a fermionic formula for
the quantity here denoted $ \langle D_{n}^{(m)}(z;q,t),e_{1}^{n}
\rangle$, and also for similar quantities in which
$D_{n}^{(m)}(z;q,t)$ is replaced by a sum over partitions $\lambda $
contained in a more general trapezoidal shape $(l^{n})+m\delta _{n}$.

Finally, we expect the analog for $m>1$ of
Conjecture~\ref{conj:nabla-En,k} to be the following.

\begin{conj}\label{conj:m-nabla-En,k}
For $1\le k \le n$,
\begin{equation}\label{e:m-nabla-En,k}
\nabla^{m} E_{n,k} =
\sum _{\substack{\lambda \subseteq m \delta _{n} \\
|\{i:\lambda _{i} = m(n-i) \}| = k}} t^{|m\delta
_{n}/\lambda |}D_n^{(m),\lambda}(z;q).
\end{equation}
\end{conj}

\subsection{Proof of Theorem~\ref{thm:m-symmetry}}
\label{s:m-symmetry}

In this section,  d-inversions and $<_{d}$ are  always defined with
respect to the given integer $m$.  We begin with a lemma.
	
\begin{lemma}\label{lemma:dinv-diff}
Let $u,v$ be cells in different columns of a tableau $T$, with $v$ to
the left of $u$.

(a) If $1-m\leq d_m(u)-d_m(v)\leq 0$, then entries $T(u)<T(v)$
contribute one more d-inversion than entries $T(u)>T(v)$.

(b) If $1\leq d_m(u)-d_m(v) \leq m$, then entries $T(u)>T(v)$ contribute
one more d-inversion than entries $T(u)<T(v)$.

(c) Otherwise, the number of d-inversions contributed by entries
$T(u)$, $T(v)$ is zero in any case.
\end{lemma}

\begin{proof}
Referring to the rule in \S \ref{s:extended-conjecture} for the number
of d-inversions contributed, we see that case (i) occurs when
$T(u)<T(v)$, and contributes $\max (0,m-|d_m(v)-d_m(u)|) = \max
(0,m-|d_m(u)-d_m(v)|)$ d-inversions.  Otherwise, case (ii) occurs and
contributes $\max (0,m-|d_m(u)-d_m(v)-1|)$ d-inversions.

If $d_m(u)>d_m(v)$ then both $d_m(u)-d_m(v)$ and $d_m(u)-d_m(v)-1$ are
non-negative.  Then case (ii) contributes one more d-inversion than
case (i), unless $d_m(u)-d_m(v)>m$, in which event both cases contribute
zero d-inversions.  If $d_m(u)\leq d_m(v)$ then both $d_m(u)-d_m(v)$ and
$d_m(u)-d_m(v)-1$ are non-positive.  Then case (i) contributes one more
d-inversion than case (ii), unless $d_m(u)-d_m(v)\leq -m$, in which event
there are again zero d-inversions in either case.
\end{proof}

Using this lemma, we can simplify the rule for counting d-inversions,
at the price of adding an overall constant.  Let $T$ be semistandard
of shape $(\lambda +(1^{n}))/\lambda $.  We say that entries $T(x) =
a$, $T(y) = b$, with $a<b$ and $x = (i,j)$, $y = (i',j')$ contribute a
{\it reduced d-inversion} if either
\begin{itemize}
\item [(i)${}'$ ] $0\leq d_{m}(y)-d_{m}(x)\leq m-1$ and $j>j'$, or
\item [(ii)${}'$ ] $1\leq d_{m}(y)-d_{m}(x)\leq m$ and $j<j'$.
\end{itemize}
Pairs of equal (positive) entries do not contribute any reduced
d-inversions.  Note that this agrees with the rule that equal entries
count as if they were unequal entries $a<b$ with $x<_{d}y$, since both
(i)${}'$  and (ii)${}'$  imply $x>_{d}y$.

Let $\dinv '_{m}(T)$ denote the number of reduced d-inversions in $T$.
Then we can rephrase Lemma~\ref{lemma:dinv-diff} as follows.

\begin{cor}\label{cor:reduced-dinv}
There is a constant $e(\nu )$ depending only on $\nu = (\lambda
+(1^{n}))/\lambda $ such that $\dinv _{m}(T) = e(\nu ) + \dinv
'_{m}(T)$ for all $T\in \SSYT (\nu )$.  Consequently,
\begin{equation}\label{e:reduced-Dnm}
D_{n}^{(m),\lambda }(z;q) = q^{e(\nu )} \sum _{T\in \SSYT (\lambda
+(1^{n})/\lambda )} q^{\dinv '_{m} (T)} z^{T}.
\end{equation}
\end{cor}

\begin{remark}
In effect, the constant $e(\nu )$ is $\dinv _{m} (T_{0})$, where
$T_{0}$ is the tableau of shape $\nu $ whose entries are all $1$
(however, $T_{0}$ is generally not a legal semistandard tableau).
\end{remark}

To prove Theorem~\ref{thm:m-symmetry}, it suffices to identify the sum
on the right-hand side of \eqref{e:reduced-Dnm} with an expression of
the form \eqref{e:Gmu-by-inv}, in this case with $n$ replaced by
$mn+1$, so $\boldmu = (\mu ^{(0)},\ldots,\mu ^{(mn)})$, and each $\mu
^{(j)}$ a partition shape in order to deduce positivity along with
symmetry.

To this end, write $\lambda = (0^{\alpha _{0}},1^{\alpha
_{1}},\ldots,mn^{\alpha _{mn}})$ with $\alpha _{0}$ defined so that
$\sum _{j}\alpha _{j} = n$.  Let $\beta $ be the permutation of
$\{0,1,\ldots,mn \}$ such that
\begin{equation}\label{e:beta}
\beta (j)\equiv -nj \pmod{mn+1}.
\end{equation}
It exists because $n$ is relatively prime to $mn+1$.  Note that for
$m=1$, $\beta $ is the identity permutation, which is why it did not
come up earlier in the proof for that case.  Define each $\mu ^{(j)}$ to
be a single column, such that
\begin{equation}\label{e:mu-j}
\mu ^{(\beta (j))} = (1^{\alpha _{j}}).
\end{equation}
We have a natural bijection between cells $x\in (\lambda
+(1^{n}))/\lambda $ and $x'\in \boldmu $, which translates column $j$
of $(\lambda +(1^{n}))/\lambda $ onto $\mu ^{(\beta (j))}$.  This
induces a bijection of semistandard tableaux in the obvious way.

Define content offsets 
\begin{equation}\label{e:s-beta-j}
s_{\beta (j)} = -nj-(mn+1)\lambda '_{j+1}.
\end{equation}
Note that the right hand side is congruent to $\beta (j)\pmod{mn+1}$,
as it should be.  With these offsets, the adjusted content of the cell
$x'\in \boldmu $ corresponding to $x = (i,j)\in (\lambda
+(1^{n}))/\lambda $ is
\begin{equation}\label{e:c-tilde-m}
\tilde{c}(x') = -(mn+1)i-nj = -i-nd_{m}(x).
\end{equation}
For any two distinct cells $x=(i,j)$ and $y = (i',j')$, we have
$0<|i-i'|<n$.  It follows that the inequalities
\begin{equation}\label{e:c-tild-range}
0<\tilde{c}(x')-\tilde{c}(y')<mn+1
\end{equation}
hold if and only if 
\begin{equation}\label{e:d-diff-range}
0\leq d_{m}(y)-d_{m}(x)\leq m,
\end{equation}
and also $i<i'$ if $d_{m}(y)-d_{m}(x) = 0$, and $i>i'$ if
$d_{m}(y)-d_{m}(x) = m$.  Moreover, since $j=j'$ implies that $mn+1$
divides $\tilde{c}(y)-\tilde{c}(x)$, these conditions imply $j\not
=j'$.  Hence we have $i<i'\Leftrightarrow j>j'$ and
$i>i'\Leftrightarrow j<j'$.  In short, inequalities
\eqref{e:c-tild-range} hold if and only if (i)${}'$ or (ii)${}'$
holds for the cells $x$, $y$.

This shows that if $S\in \SSYT (\boldmu )$ corresponds under the
natural bijection to $T\in \SSYT (\lambda +(1^{n})/\lambda )$, then
$\inv (S)  = \dinv '_{m}(T)$.  The desired result follows.
\qed 

\section{Open problems}
\label{problems}

\begin{problem}\label{prob:main}
Prove Conjecture~\ref{conj:main} and, more generally,
Conjecture~\ref{conj:m}.
\end{problem}

\begin{problem}\label{prob:m-cat}
Prove that $\langle \nabla ^{m}e_{n},e_{n} \rangle = \langle
D_{n}^{(m)}(z;q,t),e_{n} \rangle$ for $m>1$, or what is the same,
prove the combinatorial formulae \eqref{e:m-qt-Catalan-dinv} and
\eqref{e:m-qt-Catalan-hook} for $C_{n}^{(m)}(q,t)$.  This problem
might be amenable to attack by methods similar to those used in
\cite{GaHag01a,GaHag02} for the case $m=1$.
\end{problem}

\begin{problem}\label{prob:qt-sym}
Prove that $D_{n}(z;q,t) = D_{n}(z;t,q)$, and similarly for
$D_{n}^{(m)}(z;q,t)$.  This would remain a relevant combinatorial
problem even if the conjectures were proven.  By way of illustration,
although the combinatorial formula \eqref{e:qt-Catalan-theorem} for
$C_{n}(q,t)$ has been proven, no combinatorial interpretation of the
symmetry $C_{n}(q,t) = C_{n}(t,q)$ is known at present.
\end{problem}

\begin{problem}\label{prob:t=1}
Prove that $D_{n}^{(m)}(z;q,1) = D_{n}^{(m)}(z;1,q)$.  This
specialization of Problem~\ref{prob:qt-sym} should be easier than the
full $q,t$ symmetry.  Combinatorial interpretations are known
\cite{EgHaKiKr03, Loehr03a, Loehr03, LoeRem03} for the special cases
corresponding to the identities $C_{n}^{(m)}(q,1) = C_{n}^{(m)}(1,q)$,
$\langle D_{n}(z;q,1),e_{n-d}h_{d} \rangle = \langle
D_{n}(z;1,q),e_{n-d}h_{d} \rangle$ and $\langle
D_{n}^{(m)}(z;q,1),e_{1}^{n} \rangle = \langle
D_{n}^{(m)}(z;1,q),e_{1}^{n} \rangle$.
\end{problem}

\begin{problem}\label{prob:t=1/q}
Prove that 
\begin{equation}\label{e:t=1/q}
D_{n}^{(m)}(z;q,q^{-1}) = q^{-m \binom{n}{2}}
\frac{e_{n}[Z[mn+1]_{q}]}{[mn+1]_{q}}.
\end{equation}
The right-hand side of \eqref{e:t=1/q} is equal to $\nabla ^{m}
e_{n}(z)|_{t=q^{-1}}$ by \cite[Thm.~5.1]{GaHa96}.  The Catalan
specialization 
\begin{equation}\label{e:cat-t=1/q}
\langle D_{n}^{(m)}(z;q,q^{-1}), e_{n} \rangle = C_{n}^{(m)}(q,q^{-1})
\end{equation}
has been shown in \cite{Loehr03a, Loehr03}.

We remark that by the Cauchy formula,
\begin{equation}\label{e:expand-e}
\langle e_{n}[Z[mn+1]_{q}],h_{\mu } \rangle = q^{n(\mu ')}\prod _{i}
\qbinom{mn+1}{\mu _{i}}_{q} = \sum _{\substack{\lambda \subseteq
((mn+1)^{n}) \\ T\in \SSYT (\lambda +(1^{n})/\lambda, \, \mu )}}
q^{|\lambda |}.
\end{equation}
Then \eqref{e:t=1/q}  asserts that this last expression is equal to
\begin{equation}\label{e:t=1/q-detail}
[mn+1]_{q}\; \cdot  \sum _{\substack{\lambda \subseteq m\delta _{n} \\
T\in \SSYT (\lambda +(1^{n})/\lambda , \, \mu ) }} q^{|\lambda |+\dinv
_{m}(T)}.
\end{equation}
\end{problem}



\providecommand{\bysame}{\leavevmode\hbox to3em{\hrulefill}\thinspace}

\end{document}